\documentclass[authoryear,final,1p,times]{elsarticle}
\usepackage{amsmath}    
\usepackage{amssymb}    
\usepackage{amsthm}     
\usepackage{amsfonts}
\usepackage[latin1]{inputenc}
\usepackage{latexsym}
\usepackage[english]{babel}
\usepackage[T1]{fontenc}
\usepackage{graphicx}
\usepackage{xspace}
\usepackage{multirow}
\journal{Arxiv}

\begin{document}

\newcommand{\pd}    {\backslash}
\newcommand{\Pn}    {\mathbb{P}}
\newcommand{\PP}    {{\bf P}}
\newcommand{\e}     {\ensuremath{\varepsilon}}
\newcommand{\dr}    {\ensuremath{\partial }}
\newcommand{\indep} {\ensuremath{\!\!\perp\!\!\!\perp\!\!}}
\newcommand{\transp}[1] {#1 ^{\top} }
\newcommand{\C}{\mathbb{C}}
\newcommand{\N}{\mathbb{N}}
\newcommand{\R}{\mathbb{R}}
\newcommand{\Z}{\mathbb{Z}}
\newcommand{\E}{\mathbb{E}}
\newcommand{\HH}{\mathbb{H}}
\newcommand{\K}{\mathbb{K}}
\newcommand{\T}{\mathbb{T}}
\newcommand{\cA}{{\mathcal A}}
\newcommand{\cB}{{\mathcal B}}
\newcommand{\cC}{{\mathcal C}}
\newcommand{\cD}{{\mathcal D}}
\newcommand{\cE}{{\mathcal E}}
\newcommand{\cF}{{\mathcal F}}
\newcommand{\cG}{{\mathcal G}}
\newcommand{\cH}{{\mathcal H}}
\newcommand{\cI}{{\mathcal I}}
\newcommand{\cJ}{{\mathcal J}}
\newcommand{\cK}{{\mathcal K}}
\newcommand{\cL}{{\mathcal L}}
\newcommand{\cM}{{\mathcal M}}
\newcommand{\cN}{{\mathcal N}}
\newcommand{\cO}{{\mathcal O}}
\newcommand{\cP}{{\mathcal P}}
\newcommand{\cQ}{{\mathcal Q}}
\newcommand{\cR}{{\mathcal R}}
\newcommand{\cS}{{\mathcal S}}
\newcommand{\cT}{{\mathcal T}}
\newcommand{\cV}{{\mathcal V}}
\newcommand{\cU}{{\mathcal U}}
\newcommand{\cX}{{\mathcal X}}
\newcommand{\cY}{{\mathcal Y}}
\newcommand{\cZ}{{\mathcal Z}}
\def\bibname{References}
\def\refname{References}
\newcommand{\tend}  {\rightarrow}
\newcommand{\cvL}   {\stackrel{\it\mathcal Law}{\rightarrow} }
\newcommand{\cvP}   {\stackrel{\it\mathbb P}{\rightarrow} }
\newcommand{\cvps}  {\stackrel{\bf a.s.}{\rightarrow} }
\newtheorem{definition} {\bf Definition}[section]
\newtheorem{theoreme}   {\bf Theorem}[section]
\newtheorem{lemme}      {\bf lemme}[section]
\newtheorem{corollaire} {\bf Corollary}[section]
\newtheorem{proposition}{\bf Proposition}[section]
\newtheorem{PD}         {\bf Definition and Property}[section]
\newtheorem{PROPIE}     {\bf Property}[section]
\newtheorem{exemple}    {$\mathcal Example$}[section]
\newtheorem{remarque}   {$\mathcal Remark$ }[section]
\newtheorem{simulation}    {$\mathcal Simulation$}[section]
\newenvironment{preuve}        {\noindent {$\mathcal Proof\ :     $ \\}}{\hfill$\Box$\\}
\newenvironment{notation}      {\noindent {$\mathcal Notation\ :$ }\it}{}
\numberwithin{equation}{section}

\begin{frontmatter}
\title{Goodness-of-fit Tests For Elliptical And Independent Copulas Through Projection Pursuit}
\author{Jacques Touboul}
\address{Université Pierre et Marie Curie\\
Laboratoire de Statistique Théorique et Appliquée\\
jack\_touboul@hotmail.com
}
\begin{abstract}
Two goodness-of-fit tests for copulas are being investigated.
The first one deals with the case of elliptical copulas and the second one deals with independent copulas. These tests result from the expansion of the projection pursuit methodology we will introduce in the present article. This method enables us to determine on which axis system these copulas lie as well as the exact value of these very copulas in the basis formed by the axes previously determined irrespective of their value in their canonical basis.
Simulations are also presented as well as an application to real datasets.
\end{abstract}
\end{frontmatter}
\section*{ Outline of the article}
The need to describe the dependency between two or more random variables triggered the concept of copulas. Let us consider a joint cumulative distribution function (cdf) $F$ on $\R^d$ and let us consider its cdf margins $F_1$, $F_2$, ...,$F_d$, then a copula $C$ is a function such that $$F=C(F_1, F_2, ...,F_d).$$
\citet{MR0125600} is the first to have established the bases of this new theory. Several parametric families of copulas have since been defined, namely elliptical, archimedean, periodic copulas etc - see \citet{MR1462613} and \citet{MR2197664} as well as appendix \ref{CopFam} for an overview of these families.\\
Finding criterias to determine the best copula for a given problem can only be achieved through a goodness-of-fit (GOF) approach.\\
Several GOF copula approaches have so far been proposed in the literature, e.g. \citet{MR1281214}, \citet{MR2154005}, \citet{Ferm05}, \citet{Gen06}, \citet{MR2488912}, \citet{Gen09}, \citet{Mesf09}, \citet{Gen09-2}, \citet{Berg09}, \citet{MR2575418}, among others.
However, the field is still at an embryonic stage which explains the current shortage in recommendations.
In univariate distributions, the GOF assessment can be performed using for instance the well-known Kolmogorov test.
In the multivariate field, there are fewer alternatives. A simple way to build GOF approaches for multivariate random variables is to consider multi-dimensional chi-square approaches, as in for example \citet{MR2279681}.
However, these approaches present feasibility issues for high dimensional problems due to the curse of dimensionality.
In order to solve this, we will now introduce the theory of projection pursuit.

The objective of projection pursuit is to generate one or several projections providing as much information as possible about the structure of the dataset regardless of its size.\\
Once a structure has been isolated, the corresponding data are transformed through a Gaussianization. Through a recursive approach, this process is iterated to find another structure in the remaining data, until no futher structure can be evidenced in the data left at the end.\\
\citet{Frie84} and \citet{MR790553} count among the first authors who introduced this type of approaches for evidencing structures. They each describe, with many examples, how to evidence such a structure and consequently how to estimate the density of such data through two different methodologies each. Their work is based on maximizing Kullback-Leibler divergence.\\
In the present article, we will introduce a new projection pursuit methodology based on the minimisation of any $\phi$-divergence greater than the $L^1$- distance ($\phi$-PP). As we will develop later on, this way of implementing this methodology encompasses all other previous methods. This algorithm also presents the extra advantage of being more robust and more rapid from a numerical standpoint. Finally, this process allows not only to carry out GOF tests for elliptical and independent copulas but also to determine the axis system upon which these very copulas are based. 
It will also enable us to derive the exact expression of these copulas in the basis constituted by these axes.

This paper is organised as follows : section \ref{Reminders} contains preliminary definitions and properties. In section \ref{TheAlgo}, we present in details the $\phi$-projection pursuit algorithm. In section \ref{results}, we present our first results. In section \ref{GOF}, we introduce our tests. In section \ref{Simul400}, we provide two simulations pertaining to the two major situations described herein and we will study a real case.
\section{Basic theory}\label{Reminders}
\subsection{An introduction to copulas}\label{DefCop}
In this section, we will introduce the concept of copula. We will also define the family of elliptical copulas through a brief reminder of elliptical distributions - see appendix \ref{CopFam} for an overview of other families.

\subsubsection{\it  Sklar's theorem}

First, let us define a copula in $\R^d$
\begin{definition}
A $d$-dimensional copula is a joint cumulative distribution function $C$ defined on $[0,1]^d$, with uniform  margins.
\end{definition}
Moreover, the following theorem explains in what extent a copula does describe the  dependency between two or more random variables.
\begin{theoreme}[\citet{MR0125600}]Let $F$ be a joint multivariate distribution with margins $F_1$,..., $F_d$, then, there exists a copula $C$ such that
\begin{equation}
F(x_1,...,x_d)=C(F_1(x_1),...,F_d(x_d)).
\end{equation}
Moreover, if marginal cumulative distributions are continuous, then the copula is unique. Otherwise, the copula is unique on the range of values of the marginal cumulative distributions.
\end{theoreme}
\begin{remarque}
First, for any copula $C$ and any $u_i$ in $[0,1]$, $1\leq i\leq d$, we have\\
$ W(u_1,\ldots,u_d) = \max\left\{1-d+\sum\limits_{i=1}^d {u_i} , 0 \right\} \leq C(u_1,\ldots,u_d)\le \min_{j \in \{1,\ldots,d\}} u_j = M(u_1,\ldots,u_d),$\\
where $W$ and $M$ are called the Frechet-Hoeffding copula boundaries and are also copulas.\\
Moreover, we define the independent copula $\Pi$ as $\Pi(u_1,\ldots,u_d)=\Pi_{i=1}^du_i$, for any $u_i$ in $[0,1]$, $1\leq i\leq d$.
\end{remarque}
Finally, we define the density of a copula as the density associated with the cdf $C$, that we will name $c$:
\begin{definition}
Should it exist, the density of $C$ is defined by
$c(u_1,...,u_d)=\frac{\dr^d}{\dr u_1...\dr u_d}C(u_1,...,u_d)$, for any $u_i$ in $[0,1]$, $1\leq i\leq d$.
\end{definition}
\subsubsection{\it  The Gaussian copula}

The Gaussian copula can be used in several fields. For example, many credit models are built from this copula, which also presents the property to make extreme values (minimal or maximal) independent - in the limit ; see \citet{MR1462613} for more details.
For example, in $\R^2$, it is derived from the bivariate normal distribution and from Sklar's theorem. Defining $\Psi_\rho$ as the standard bivariate normal cumulative distribution function with $\rho$ correlation, the Gaussian copula function is
$C_\rho(u,v) = \Psi_\rho \left(\Psi^{-1}(u), \Psi^{-1}(v) \right)$
where $u, v \in [0,1]$ and where $\Psi$ is the standard normal cumulative distribution function.
Then, the copula density function is :
$$c_\rho(u,v) = \frac{\psi_{X,Y, \rho} (\Psi^{-1}(u), \Psi^{-1}(v) )}
{\psi(\Psi^{-1}(u)) \psi(\Psi^{-1}(v))}$$ where 
$\psi_{X,Y, \rho}(x,y) = \frac{1}{2 \pi\sqrt{1-\rho^2}} \exp \left ( -\frac{1}{2(1-\rho^2)}  \left [{x^2+y^2} -2\rho xy  \right ] \right )$ is the density function for the standard bivariate Gaussian with pearson product-moment correlation coefficient  $\rho$ and where $\psi$ is the standard normal density. This definition can obviously be extended to $\R^d$.

\subsubsection{\it  The elliptical copula}

Let us begin with defining the class of elliptical distributions and its properties - see also \citet{MR0629795}, \citet{MR2061237} :
\begin{definition}
$X$ is said to abide by a multivariate elliptical distribution, denoted $X\sim E_d(\mu,\Sigma,\xi_d)$, if $X$ has the following density, for any $x$ in $\R^d$ :

$\ \ \ \ \ \ \ \ \ \ \ \ \ \ \ \ \ \ \ \ \ \ \ \ \ \ \ \ $ $ f_X(x)=\frac{\alpha_d}{|\Sigma|^{1/2}}\xi_d\Big(\frac{1}{2}(x-\mu)'\Sigma^{-1}(x-\mu)\Big)$\\
$\bullet$ where $\Sigma$ is a $d\times d$ positive-definite matrix and where $\mu$ is a $ d$-column vector,\\
$\bullet$ where $\xi_d$ is referred as the "density generator",\\
$\bullet$ where $\alpha_d$ is a normalisation constant, such that $\alpha_d=\frac{\Gamma(d/2)}{(2\pi)^{d/2}}\Big(\int_0^\infty x^{d/2-1}\xi_d(x)dx\Big)^{-1}$,\\ with $\int_0^\infty x^{d/2-1}\xi_d(x)dx<\infty$.
\end{definition}
\begin{PROPIE}\label{ElliptProp}
1/ For any $X\sim E_d(\mu,\Sigma,\xi_d)$, for any $m\times d$ matrix with rank $m\leq d$, $A$, and for any $m$-dimensional vector $b$, we have $AX+b\sim E_m(A\mu+b,A\Sigma A',\xi_m)$.\\
Therefore, any marginal density of multivariate elliptical distribution is elliptical, i.e. \\
$X=(X_1,X_2,...,X_d)\sim E_d(\mu,\Sigma,\xi_d)\ \Rightarrow\ X_i\sim E_1(\mu_i,\sigma^2_i,\xi_1),$ $1\leq i\leq d$, with $f_{X_i}(x)= \frac{\alpha_1}{\sigma_i}\xi_1\Big(\frac{1}{2}(\frac{x-\mu_i}{\sigma})^2\Big)$.\\
2/ Corollary 5 of \citet{MR0629795} states that conditional densities with elliptical distributions are also elliptical. Indeed, if $X=(X_1,X_2)'\sim E_d(\mu,\Sigma,\xi_d)$, with $X_1$ (resp. $X_2$) of size $d_1<d$ (resp. $d_2<d$), then $X_1/(X_2=a)\sim E_{d_1}(\mu',\Sigma',\xi_{d_1})$ with $\mu'=\mu_1+\Sigma_{12}\Sigma_{22}^{-1}(a-\mu_2)$ and $\Sigma'=\Sigma_{11}-\Sigma_{12}\Sigma_{22}^{-1}\Sigma_{21},$
with $\mu=(\mu_1,\mu_2)$ and $\Sigma=(\Sigma_{ij})_{1\leq i,j\leq 2}$.
\end{PROPIE}
\begin{remarque}\label{implyEstimBounded}
\citet{MR2061237} shows that multivariate Gaussian distributions derive from $\xi_d(x)=e^{-x}$. They also show that if $X=(X_1,...,X_d)$ has an elliptical density such that its marginals verify $E(X_i)<\infty$ and $E(X_i^2)<\infty$ for $1\leq i\leq d,$ then $\mu$ is the mean of $X$ and $\Sigma$ is a multiple of the covariance matrix of $X$. Consequently, from now on, we will assume this is indeed the case.
\end{remarque}
\begin{definition}
Let $t$ be an elliptical density on $\R^k$ and let $q$ be an elliptical density on $\R^{k'}$.
The elliptical densities $t$ and $q$ are said to belong to the same family of elliptical densities, if their generating densities are $\xi_k$ and $\xi_{k'}$ respectively, which belong to a common given family of densities.
\end{definition}
\begin{exemple}
Consider two Gaussian densities  $\cN(0,1)$ and $\cN((0,0),Id_2)$. They are said to belong to the same elliptical family as they both present  $x\mapsto e^{-x}$ as generating density.
\end{exemple}
Finally, let us introduce the definition of an elliptical copula which generalizes the above overview of the Gaussian copula :
\begin{definition}
Elliptical copulas are the copulas of elliptical distributions.
\end{definition}
\subsection{Brief introduction to the $\phi$-projection pursuit methodology ($\phi$-PP)}\label{OurProposal}
Let us first introduce the concept of $\phi-$divergence.
\subsubsection{\it  The concept of $\phi-$divergence}

Let $\varphi$ be a strictly convex function defined by $\varphi: \overline{\R^+}\to\overline{\R^+},$ and such that $\varphi(1)=0$. We define a $\phi-$divergence of $P$ from $Q$ - where $P$ and $Q$ are two probability distributions over a space $\Omega$ such that $Q$ is absolutely continuous with respect to $P$ - by  $${D_{\phi}}(Q,P)=\int\varphi(\frac{dQ}{dP})dP$$
or ${D_{\phi}}(q,p)=\int\varphi(\frac{q(x)}{p(x)})p(x)dx$, if $P$ and $Q$ present $p$ and $q$ as density respectively.\\
Throughout this article, we will also assume that $\varphi(0)<\infty$, that $\varphi'$ is continuous and that this divergence is greater than the $L^1$ distance - see also Appendix \ref{FiDiv} page \pageref{FiDiv}.

\subsubsection{\it  Functioning of the algorithm}\label{MatrixAlgo}

Let $f$ be a density on $\R^d$. We define an instrumental density $g$ with the same mean and variance as  $f$.
We start with performing the ${D_{\phi}}(g,f)=0$ test; should this test turn out to be positive, then $f=g$ and the algorithm stops, otherwise, the first step of our algorithm consists in defining a vector $a_1$ and a density $g^{(1)}$  by 
\begin{equation}\label{PhiDefSequMethod}
a_1\ =\ arg\inf_{a\in\R^d_*}\ {D_{\phi}}(g\frac{f_a}{g_a},f)\text{ and }g^{(1)}=g\frac{f_{a_1}}{g_{a_1}}
\end{equation}
where $\R^d_*$ is the set of non null vectors of $\R^d$ and $f_a$ (resp. $g_a$) stands for the density of $\transp aX$ (resp. $\transp aY$) when $f$ (resp. $g$) is the density of $X$ (resp. $Y$).\\
In our second step, we will replace $g$ with $g^{(1)}$ and we will repeat the first step. \\And so on, by iterating this process, we will end up obtaining a sequence  $(a_1,a_2,...)$ of vectors in  $\R^d_*$  and a sequence of densities  $g^{(i)}$.
We will thus prove that the underlying structures of $f$ evidenced through this method are identical to the ones obtained through projection pursuit methodologies based on Kullback-Leibler divergence maximisation, such as Huber's method - see appendix \ref{a1solve4}. We will also evidence the above structures, which will enable us to infer more information on $f$ - see example below.
\begin{remarque}
First, to obtain an approximation of $f$, we stop our algorithm when the divergence equals zero, i.e. we stop when ${D_{\phi}}(g^{(j)},f)=0$ since it implies  $g^{(j)}=f$ with $j\leq d$, or when our algorithm reaches the $d^{th}$ iteration, i.e. we approximate $f$ with $g^{(d)}$.\\
Second, we get ${D_{\phi}}(g^{(0)},f)\geq {D_{\phi}}(g^{(1)},f)\geq.....\geq 0$ with $g^{(0)}=g$.\\
Finally, the specific form of the relationship (\ref{PhiDefSequMethod}) establishes that we deal with M-estimation. We can therefore state that our method is more robust than projection pursuit methodologies based on Kullback-Leibler divergence maximisation - see \citet{BAR-YOHAI}, \citet{TOMA} as well as \citet{RobStat}.
\end{remarque}
\noindent At present, let us study the following example:
\begin{exemple}Let $f$ be a density defined on $\R^3$ by $f(x_1,x_2,x_3)=n(x_1,x_2)h(x_3)$, with $n$ being a bi-dimensional Gaussian density, and $h$ being a non Gaussian density. Let us also consider $g$, a Gaussian density with the same mean and variance as $f$.\\
Since $g(x_1,x_2/x_3)=n(x_1,x_2)$, we have ${D_{\phi}}(g\frac{f_3}{g_3},f)={D_{\phi}}(n.f_3,f)={D_{\phi}}(f,f)=0$ as $f_3=h$, i.e. the function $a\mapsto {D_{\phi}}(g\frac{f_a}{g_a},f)$ reaches zero for $e_3=(0,0,1)'$ - where $f_3$ and $g_3$ are the third marginal densities of $f$ and $g$ respectively. We therefore obtain $g(x_1,x_2/x_3)=f(x_1,x_2/x_3)$.
\end{exemple}
To recapitulate our method, if ${D_{\phi}}(g,f)=0$, we derive $f$ from the relationship $f=g$; should a sequence $(a_i)_{i=1,...j}$, $j<d$, of vectors in $\R^d_*$ defining $g^{(j)}$ and such that ${D_{\phi}}(g^{(j)},f)=0$ exist, then $f(./\transp{a_i}x, 1\leq i\leq j)=g(./\transp{a_i}x, 1\leq i\leq j)$, i.e. $f$ coincides with $g$ on the complement  of the vector subspace generated by the family $\{a_i\}_{i=1,...,j}$ - see also section \ref{TheAlgo} for a more detailed explanation. 

In the remaining of the study of the algorithm, after having clarified the choice of $g$,  we will consider the statistical solution to the representation problem, assuming that $f$ is unknown and that $X_1$, $X_2$,... $X_m$ are i.i.d. with density $f$. We will provide asymptotic results pertaining to the family of optimizing vectors $a_{k,m}$ - that we will define more precisely below - as $m$ goes to infinity.
Our results also prove that the empirical representation scheme converges towards the theoretical one.

\section{ The algorithm}\label{TheAlgo}
\subsection{ The model}\label{modelSection}

Let $f$ be a density on $\R^d$. We assume there exists $d$ non null linearly independent vectors $a_j$, with $1\leq j\leq d,$ of $\R^d$, such that
\begin{equation}\label{f-def}
f(x)=n(\transp{a_{j+1}}x,...,\transp{a_{d}}x)h(\transp{a_{1}}x,...,\transp{a_{j}}x)
\end{equation}
with $j<d$, $n$ being an elliptical density on $\R^{d-j}$ and with $h$ being a density on $\R^{j}$, which does not belong to the same family as $n$. Let $X=(X_{1},...,X_{d})$ be a vector with $f$ as density.\\
We define $g$ as an elliptical distribution with the same mean and variance as $f$.\\
For simplicity, let us assume that the family $\{a_j\}_{1\leq j\leq d}$ is the canonical basis of $\R^d$:\\
The very definition of $f$ implies that $(X_{j+1},...,X_{d})$ is independent from $(X_{1},...,X_{j})$. Hence,  the density of $(X_{j+1},...,X_{d})$ given $(X_{1},...,X_{j})$ is $n$.\\
Let us assume that ${D_{\phi}}(g^{(j)},f)=0,$ for some $j\leq d$. We then get $\frac{f(x)}{f_{a_1}f_{a_2}...f_{a_j}}=\frac{g(x)}{g^{(1-1)}_{a_1}g^{(2-1)}_{a_2}...g^{(j-1)}_{a_j}}$, since,  by induction, we have $g^{(j)}(x)=g(x)\frac{f_{a_1}}{g^{(1-1)}_{a_1}}\frac{f_{a_2}}{g^{(2-1)}_{a_2}}...\frac{f_{a_j}}{g^{(j-1)}_{a_j}}$.\\
Consequently, through lemma \ref{TrucBidule} and through the fact that the conditional densities with elliptical distributions are also elliptical, as well as through the above relationship, we can infer that

$\ \ \ \ \ \ \  \ \ \ \ \ \ \  \ \ $    $n(\transp{a_{j+1}}x,.,\transp{a_{d}}x)=f(./\transp{a_i}x, 1\leq i\leq j)=g(./\transp{a_i}x, 1\leq i\leq j).$\\ In other words, $f$ coincides with $g$ on the complement  of the vector subspace generated by the family $\{a_i\}_{i=1,...,j}$.

Now, if the family $\{a_j\}_{1\leq j\leq d}$ is no longer the canonical basis of $\R^d$, then this family is again a basis of $\R^d$. Hence, lemma \ref{ChangBasis} implies that
\begin{equation}\label{RelElli36}
g(./\transp{a_{1}}x,...,\transp{a_{j}}x)=n(\transp{a_{j+1}}x,...,\transp{a_{d}}x)=f(./\transp{a_{1}}x,...,\transp{a_{j}}x)
\end{equation}
which is equivalent to ${D_{\phi}}(g^{(j)},f)=0$, since by induction  $g^{(j)}=g\frac{f_{a_1}}{g^{(1-1)}_{a_1}}\frac{f_{a_2}}{g^{(2-1)}_{a_2}}...\frac{f_{a_j}}{g^{(j-1)}_{a_j}}$.\\
The end of our algorithm implies that $f$ coincides with $g$ on the complement  of the vector subspace generated by the family $\{a_i\}_{i=1,...,j}$.
Therefore, the nullity of the $\phi-$divergence provides us with information on the density structure.\\
\noindent In summary, the following proposition clarifies our choice of $g$ which depends on the family of distribution one wants to find in $f$ :
\begin{proposition}\label{Pb1}With the above notations, ${D_{\phi}}(g^{(j)},f)=0$ is equivalent to
$$
g(./\transp{a_{1}}x,...,\transp{a_{j}}x)=f(./\transp{a_{1}}x,...,\transp{a_{j}}x)
$$
\end{proposition}
More generally, the above proposition leads us to defining the co-support of $f$ as the vector space generated by the vectors $a_{1},...,a_{j}$.
\begin{definition}
Let $f$ be a density on $\R^d$. We define the co-vectors of $f$ as the sequence of vectors $a_{1},...,a_{j}$ which solves the problem ${D_{\phi}}(g^{(j)},f)=0$ where $g$ is an elliptical distribution with the same mean and variance as $f$.
We define the co-support of $f$ as the vector space generated by the vectors $a_{1},...,a_{j}$.
\end{definition}
\begin{remarque}Any $(a_i)$ family defining $f$ as in (\ref{f-def}), is an orthogonal basis of $\R^d$ - see lemma \ref{Base}
\end{remarque}
\subsection{Stochastic outline of our algorithm}\label{UseSample}
Let $X_1$, $X_2$,..,$X_m$ (resp. $Y_1$, $Y_2$,..,$Y_m$) be a sequence of $m$ independent random vectors with the same density $f$ (resp. $g$). 
As customary in nonparametric $\phi-$divergence optimizations, all estimates of $f$ and $f_a$, as well as all uses of Monte Carlo methods are being performed using subsamples $X_1$, $X_2$,..,$X_n$ and $Y_1$, $Y_2$,..,$Y_n$ - extracted respectively from $X_1$, $X_2$,..,$X_m$ and $Y_1$, $Y_2$,..,$Y_m$ - since the estimates are bounded below by some positive deterministic sequence $\theta_m$ - see Appendix \ref{truncSample}.\\
Let $\Pn_n$ be the empirical measure based on the subsample $X_1$, $X_2$,.,$X_n$. Let $f_n$ (resp.  $f_{a,n}$ for any $a$ in $\R^d_*$) be the kernel estimate of $f$ (resp.  $f_a$), which is built from $X_1$, $X_2$,..,$X_n$ (resp. $\transp aX_1$, $\transp aX_2$,..,$\transp aX_n$).\\
As defined in section \ref{OurProposal}, we introduce the following sequences $(a_{k})_{k\geq 1}$ and $(g^{(k)})_{k\geq 1}$:
\begin{eqnarray}\label{VraiDefOfAK}
&&\text{$\bullet$ $a_{k}$ is a non null vector of $\R^d$ such that $a_{k}=arg\min_{a\in\R^d_*} {D_{\phi}}(g^{(k-1)}\frac{f_a}{g^{(k-1)}_a},f)$}\\
&&\text{$\bullet$ $g^{(k)}$ is the density such that $g^{(k)}=g^{(k-1)}\frac{f_{a_{k}}}{g^{(k-1)}_{a_{k}}}$ with $g^{(0)}=g$}\nonumber
\end{eqnarray}
The stochastic setting up of the algorithm uses $f_n$ and $g_n^{(0)}=g$ instead of $f$ and $g^{(0)}=g$, since $g$ is known. Thus, at the first step, we build the vector $\check a_1$ which minimizes the $\phi-$divergence between $f_n$ and $g\frac{f_{a,n}}{g_{a}}$ and which estimates $a_1$.\\
Proposition \ref{QuotientDonneLoi} and lemma \ref{toattain} enable us to minimize the $\phi-$divergence between $f_n$ and $g\frac{f_{a,n}}{g_{a}}$. Defining $\check a_1$ as the argument of this minimization, proposition \ref{KernelpConv2} shows us that this vector tends to $a_1$.\\
Finally, we define the density $\check g^{(1)}_m$ as $\check g^{(1)}_m=g\frac{f_{\check a_1,m}}{g_{\check a_1}}$ which estimates $g^{(1)}$ through theorem \ref{KernelKRessultatPricipal}.\\
Now, from the second step and as defined in section \ref{OurProposal}, the density $g^{(k-1)}$ is unknown. Consequently, once again, we have to truncate the samples.\\
All estimates of $f$ and $f_a$ (resp. $g^{(1)}$ and $g_a^{(1)}$) are being performed using a subsample $X_1$, $X_2$,..,$X_n$ (resp. $Y_1^{(1)}$, $Y_2^{(1)}$,..,$Y_n^{(1)}$) extracted from $X_1$, $X_2$,..,$X_m$ (resp. $Y_1^{(1)}$, $Y_2^{(1)}$,..,$Y_m^{(1)}$ - which is a sequence of $m$ independent random vectors with the same density $g^{(1)}$) such that the estimates are bounded below by some positive deterministic sequence $\theta_m$ (see Appendix \ref{truncSample}).\\
Let $\Pn_n$ be the empirical measure based on the subsample $X_1$, $X_2$,..,$X_n$. Let $f_n$ (resp. $g_n^{(1)}$, $f_{a,n}$, $g_{a,n}^{(1)}$ for any $a$ in $\R^d_*$) be the kernel estimate of $f$ (resp. $g^{(1)}$ and $f_a$ as well as $g_a^{(1)}$) which is built from $X_1$, $X_2$,..,$X_n$ (resp. $Y_1^{(1)}$, $Y_2^{(1)}$,..,$Y_n^{(1)}$ and $\transp aX_1$, $\transp aX_2$,..,$\transp aX_n$ as well as $\transp aY_1^{(1)}$, $\transp aY^{(1)}_2$,..,$\transp aY_n^{(1)}$).
The stochastic setting up of the algorithm uses $f_n$ and $g_n^{(1)}$ instead of $f$ and $g^{(1)}$.
Thus, we build the vector $\check a_2$ which minimizes the $\phi-$divergence between $f_n$ and $g_n^{(1)}\frac{f_{a,n}}{g_{a,n}^{(1)}}$ - since $g^{(1)}$ and $g^{(1)}_a$ are unknown - and which estimates $a_2$.
Proposition \ref{QuotientDonneLoi} and lemma \ref{toattain} enable us to minimize the $\phi-$divergence between $f_n$ and $g_n^{(1)}\frac{f_{a,n}}{g_{a,n}^{(1)}}$. Defining $\check a_2$ as the argument of this minimization, proposition \ref{KernelpConv2} shows that this vector tends to $a_2$ in $n$. Finally, we define the density $\check g^{(2)}_n$ as $\check g^{(2)}_n=g_n^{(1)}\frac{f_{\check a_2,n}}{g^{(1)}_{\check a_2,n}}$ which estimates $g^{(2)}$ through theorem \ref{KernelKRessultatPricipal}.\\
And so on, we will end up obtaining a sequence  $(\check a_1,\check a_2,...)$ of vectors in  $\R^d_*$ estimating the co-vectors of $f$ and a sequence of densities $(\check g^{(k)}_n)_k$ such that $\check g^{(k)}_n$ estimates $g^{(k)}$ through theorem \ref{KernelKRessultatPricipal}.

\section{Results}\label{results}
\subsection{ Hypotheses on $f$}\label{HypoF}

Let $X_1$, $X_2$,..,$X_m$ (resp. $Y_1$, $Y_2$,..,$Y_m$) be a sequence of $m$ independent random vectors with the same density $f$ (resp. $g$). 
As customary in nonparametric $\phi-$divergence optimizations, all estimates of $f$ and $f_a$ as well as all uses of Monte Carlo methods are being performed using subsamples $X_1$, $X_2$,..,$X_n$ and $Y_1$, $Y_2$,..,$Y_n$ - extracted respectively from $X_1$, $X_2$,..,$X_m$ and $Y_1$, $Y_2$,..,$Y_m$ - since the estimates are bounded below by some positive deterministic sequence $\theta_m$ - see appendix \ref{truncSample}.\\
Let $\Pn_n$ be the empirical measure of the subsample $X_1$, $X_2$,.,$X_n$. Let $f_n$ (resp.  $f_{a,n}$ for any $a$ in $\R^d_*$) be the kernel estimate of $f$ (resp.  $f_a$), which is built from $X_1$, $X_2$,..,$X_n$ (resp. $\transp aX_1$, $\transp aX_2$,..,$\transp aX_n$).\\
At present, let us define the set of hypotheses on $f$.\\
Discussion on several of these hypotheses can be found in appendix \ref{DiscussHyp}.\\
In the remaining of this section, to be more legible we replace $g$ with $g^{(k-1)}$. Let

$\Theta =\R^d,\ \Theta^{{D_{\phi}}} =\{b\in\Theta\ |\ \ \int\varphi^*(\varphi'(\frac{g(x)}{f(x)}\frac{f_b(\transp bx)}{g_b(\transp bx)}))d{\bf P}<\infty\},$

$M(b,a,x)=\int\varphi'(\frac{g(x)}{f(x)}\frac{f_b(\transp bx)}{g_b(\transp bx)})g(x)\frac{f_a(\transp ax)}{g_a(\transp ax)}dx-\ \varphi^*(\varphi'(\frac{g(x)}{f(x)}\frac{f_b(\transp bx)}{g_b(\transp bx)})),$

$\Pn_n M(b,a)=\int M(b,a,x)d\Pn_n,$ ${\bf P} M(b,a)=\int M(b,a,x) d{\bf P},$\\
where $\PP $ is the probability measure presenting $f$ as density.\\
Similarly as in chapter $V$ of \citet{MR1652247}, let us define :\\
$(H1)$ : $\text{For all $\e>0$, there is $\eta>0$, such that for all $c\in\Theta^{D_{\phi}} $ verifying $\|c-a_k\|\geq \e,$}$

\hspace{1cm}$\text{ we have }{\bf P} M(c,a)-\eta>{\bf P} M(a_k,a),\text{ with }a\in\Theta.$\\
$(H2)$ : $\text{$\exists$ $Z<0$, $n_0>0$ such that }(n\geq n_0\ \Rightarrow \ \sup_{a\in\Theta}\sup_{c\in\{\Theta^{D_{\phi}} \}^c}\Pn_nM(c,a)<Z)$\\
$(H3)$ : $\text{There exists $V$, a neighbourhood of $a_k$, and $H$, a positive function, such that, for all $c\in V$,}$

\hspace{1cm}$\text{ we have }|M(c,a_k,x)|\leq H(x) ({\bf P} -a.s.)\text{ with }{\bf P} H<\infty,$\\
$(H4)$ : $\text{There exists $V$, a neighbourhood of $a_k$, such that for all $\e$, there exists a $\eta$ such that for all $c \in V$}$

\hspace{1cm}$\text{ and $a\in\Theta$, verifying }\|a-a_k\|\geq \e,\text{ we have }{\bf P} M(c,a_k)<{\bf P} M(c,a)-\eta.$\\
Putting $I_{a_k}=\frac{\dr^2}{\dr a^2}   {D_{\phi}} (g\frac{f_{a_k}}{g_{a_k}},f)$, let us consider now four new hypotheses:\\
$(H5)$ : ${\bf P}\|\frac{\dr}{\dr b}M(a_k,a_k)\|^2$ and ${\bf P}\|\frac{\dr}{\dr a}M(a_k,a_k)\|^2$ are finite and the expressions ${\bf P}\frac{\dr^2}{\dr b_i\dr b_j}M(a_k,a_k)$ and $I_{a_k}$

\hspace{1cm}  exist and are invertible.\\
$(H6)$ : There exists $k$ such that  ${\bf P} M(a_k,a_k)= 0$.\\
$(H7)$ : $(Var_{{\bf P}}(M(a_k,a_k)))^{1/2}$ exists and is invertible.\\
$(H0)$ : $f$ and $g$ are assumed to be positive and bounded and such that $K(g,f)\geq\int|f(x)-g(x)|dx$

\hspace{1cm} where $K$ is the Kullback-Leibler divergence.
\subsubsection{\it   Estimation of the first co-vector of $f$}\label{Estimofa1}
Let $\cR$ be the class of all positive functions $r$ defined on $\R$ and such that $g(x)r(\transp ax)$ is a density on $\R^d$ for all $a$ belonging to $\R^d_*$. The following proposition shows that there exists a vector $a$ such that $\frac{f_a}{g_a}$ minimizes ${D_{\phi}}(gr,f)$ in $r$:
\begin{proposition} \label{lemmeHuberModifprop}
There exists a vector $a$ belonging to $\R^d_*$ such that
$$
arg\min_{r\in\cR}{D_{\phi}}(gr,f)=\frac{f_a}{g_a}\text{ and }r(\transp ax)=\frac{f_a(\transp ax)}{g_a(\transp ax)}
$$
\end{proposition}
\begin{remarque}\label{criteria}
This proposition proves that $a_1$ simultaneously optimises (\ref{DefSequMethodH1}), (\ref{DefSequMethodH2}) and (\ref{PhiDefSequMethod}). In other words, it proves that the underlying structures of $f$ evidenced through our method are identical to the ones obtained through projection pursuit methodologies based on Kullback-Leibler divergence maximisation, such as Huber's methods - see appendix \ref{PP-old}.
\end{remarque}
Following \citet{MR2054155}, let us introduce the estimate of ${D_{\phi}}(g\frac{f_{a,n}}{g_{a}},f_{n})$, through 

$\ \ \ \ \ \ \ \ \ \ \ \ \ \ \ \ \ \ \ \ \ \ \ \ \ \ \ \ \ \ \ \ \ \ \ \ \ \ \ \ $ $\check   {D_{\phi}}(g\frac{f_{a,n}}{g_{a}},f_{n})=\int M(a,a,x) d\Pn_n(x)$
\begin{proposition}\label{info}
Let $\check  a$ be such that $\check  a   :=  arg\inf_{a\in\R^d_*}\check   {D_{\phi}}(g\frac{f_{a,n}}{g_{a}},f_{n}).$\\
Then, $\check  a$ is a strongly convergent estimate of $a$, as defined in proposition \ref{lemmeHuberModifprop}.
\end{proposition}
\noindent Let us also introduce the following sequences  $(\check a_{k})_{k\geq 1}$ and $(\check g^{(k)}_{n})_{k\geq 1}$, for any given $n$ - see section \ref{UseSample}.:\\
$\bullet$ $\check a_{k}$ is an estimate of $a_{k}$ as defined in proposition \ref{info}  with $\check g^{(k-1)}_{n}$ instead of $g$,\\
$\bullet$ $\check g^{(k)}_{n}$ is such that $\check g^{(0)}_{n}=g$, $\check g^{(k)}_{n}(x)=\check g^{(k-1)}_{n}(x)\frac{f_{\check a_k,n}(\transp {\check a_k}x)}{[\check g^{(k-1)}]_{\check a_k,n}(\transp {\check a_k}x)}$, i.e.
$\check g^{(k)}_{n}(x)=g(x)\Pi_{j=1}^k\frac{f_{\check a_j,n}(\transp {\check a_j}x)}{[\check g^{(j-1)}]_{\check a_j,n}(\transp {\check a_j}x)}$.\\
We also note that $\check g^{(k)}_n$ is a density.
\subsubsection{\it   Convergence study at the  $k^{\text{th}}$ step of the algorithm:}
In this paragraph, we show that the sequence $(\check a_k)_n$ converges towards $a_k$ and that the sequence $(\check g^{(k)}_n)_n$ converges towards $g^{(k)}$.\\ 
Let $\check c_n(a)    =\ arg\sup_{c\in\Theta }\ \Pn_nM(c,a),$ with $a\in\Theta$,
and $\check \gamma_n   =\ arg\inf_{a\in\Theta }\ \sup_{c\in\Theta }\ \Pn_nM(c,a)$.
We state
\begin{proposition}\label{KernelpConv2}	Both $\sup_{a\in\Theta}\|\check c_n(a)-a_k\|$  and  $\check \gamma_n $ converge toward $a_k$ a.s.
\end{proposition}
\noindent Finally, the following theorem shows that $\check  g^{(k)}_n$ converges almost everywhere towards $g^{(k)}$:
\begin{theoreme}\label{KernelKRessultatPricipal}It holds
$\check  g^{(k)}_n\to_n   g^{(k)}\ a.s.$
\end{theoreme}
\subsubsection{\it   Testing of the criteria}\label{Test}
In this paragraph, through a test of our criteria,  namely $a\mapsto {D_{\phi}}(\check g^{(k)}_n\frac{f_{a,n}}{[\check g^{(k)}]_{a,n}},f_n)$, we will build a stopping rule for this procedure. First, the next theorem enables us to derive the law of our criteria:
\begin{theoreme} \label{KernelLOIDUCRITERE} For a fixed $k$, we have 

$\sqrt n(Var_{\PP}(M(\check c_n(\check \gamma_n),\check \gamma_n)))^{-1/2}(\Pn_nM(\check c_n(\check \gamma_n),\check \gamma_n)-\Pn_nM(a_k,a_k)) \cvL \cN(0,I)$,\\ where $k$ represents the $k^{th}$ step of our algorithm  and where $I$ is the identity matrix in $\R^d$.
\end{theoreme}
Note that $k$ is fixed in theorem \ref{KernelLOIDUCRITERE} since $\check \gamma_n   =\ arg\inf_{a\in\Theta }\ \sup_{c\in\Theta }\ \Pn_nM(c,a)$ where $M$ is a known function of $k$ - see section \ref{HypoF}. Thus, in the case when ${D_{\phi}}(g^{(k-1)}\frac{f_{a_k}}{g^{(k-1)}_{a_k}},f)= 0$, we obtain
\begin{corollaire} \label{KernelLOIDUCRITERE2} 
We have $\sqrt n(Var_{\PP}(M(\check c_n(\check \gamma_n),\check \gamma_n)))^{-1/2}\Pn_nM(\check c_n(\check \gamma_n),\check \gamma_n) \cvL \cN(0,I)$.
\end{corollaire} 
\noindent Hence, we propose the test of the null hypothesis 

$\text{$(H_0)$ : ${D_{\phi}}(g^{(k-1)}\frac{f_{a_k}}{g^{(k-1)}_{a_k}},f)= 0$ versus the alternative $(H_1)$ :  ${D_{\phi}}(g^{(k-1)}\frac{f_{a_k}}{g^{(k-1)}_{a_k}},f)\not= 0$.}$\\
Based on this result, we stop the algorithm, then, defining $a_k$ as the last vector generated, we derive from corollary \ref{KernelLOIDUCRITERE2} a $\alpha$-level confidence ellipsoid around $a_k$, namely 

$\cE_k=\{b\in\R^d;\ \sqrt n(Var_{\PP}(M(b,b)))^{-1/2}\Pn_nM(b,b)\leq q_{\alpha}^{\cN(0,1)} \}$\\
where $q_{\alpha}^{\cN(0,1)}$ is the quantile of a $\alpha$-level reduced centered normal distribution and where $\Pn_n$ is the empirical measure araising from a realization of the sequences $(X_1,\ldots,X_n)$ and $(Y_1,\ldots,Y_n)$.\\ Consequently, the following corollary provides us with a confidence region for the above test:
\begin{corollaire}\label{KernelLOIDUCRITERE2coro}
$\cE_k$ is a confidence region for the test of the null hypothesis $(H_0)$ versus $(H_1)$.
\end{corollaire}
\section{Goodness-of-fit tests}\label{GOF}
\subsection{The basic idea}\label{ExmapleCopula}
Let $f$ be a density defined on $\R^2$. Let us also consider $g$, a known elliptical density with the same mean and variance as $f$.
Let us also assume that the family $(a_i)$ is the canonical basis of $\R^2$ and that $D_{\phi}(g^{(2)},f)=0$.\\
Hence, since lemma \ref{TrucBidule} page \pageref{TrucBidule} implies that $g^{(j-1)}_{a_j}=g_{a_j}$ if $j\leq d$, we then have $g^{(2)}(x)=g(x)\frac{f_1}{g_1}\frac{f_2}{g^{(1)}_2}=g(x)\frac{f_1}{g_1}\frac{f_2}{g_2}$.
Moreover, we get $f$ with $g^{(2)}=f$, as derived from property \ref{Phimini} page \pageref{Phimini}.\\
Consequently, $f=g(x)\frac{f_1}{g_1}\frac{f_2}{g_2},$ i.e. $\frac{f}{f_1f_2}=\frac{g}{g_1g_2}$, and then $$\frac{\dr^2}{\dr x\dr y}C_f=\frac{\dr^2}{\dr x\dr y}C_g$$ where $C_f$ (resp. $C_g$) is the copula of $f$ (resp. $g$). \\
More generally, if $f$ is defined on $\R^d$, then the family $(a_i)$ is once again free - see lemma \ref{imFree} page \pageref{imFree} -, i.e.  the family $(a_i)$ is once again a basis of $\R^d$. The relationship ${D_{\phi}}(g^{(d)},f)=0$ therefore implies that $g^{(d)}=f$, i.e. for any $x\in\R^d$,
$f(x)=g^{(d)}(x)=g(x)\Pi_{k=1}^d\frac{f_{a_k}(\transp {a_k}x)}{[g^{\{k-1\}}]_{a_k}(\transp {a_k}x)}=g(x)\Pi_{k=1}^d\frac{f_{a_k}(\transp {a_k}x)}{g_{a_k}(\transp {a_k}x)}$ since lemma \ref{TrucBidule} page \pageref{TrucBidule} implies that $g^{(k-1)}_{a_k}=g_{a_k}$ if $k\leq d$. In other words, for any $x\in\R^d$, it holds
\begin{equation}
\frac{g(x)}{\Pi_{k=1}^d g_{a_k}(\transp {a_k}x)}=\frac{f(x)}{\Pi_{k=1}^df_{a_k}(\transp {a_k}x)}\label{Copule-1}
\end{equation}
Finally, putting $A=(a_1,...,a_d)$ and defining vector $y$ (resp. density $\tilde f$, copula $\tilde C_f$ of $\tilde f$, density $\tilde g$, copula $\tilde C_g$ of $\tilde g$) as the expression of vector $x$ (resp. density $f$, copula $C_f$ of $ f$, density $g$, copula $C_g$ of $g$) in basis $A$, 
then, the following proposition provides us with the density associated with the copula of $f$ as being equal to the density associated with the copula of $g$ in basis $A$ :
\begin{proposition}\label{PropCop}
With the above notations, should a sequence $(a_i)_{i=1,...d}$ of not null vectors in $\R^d_*$ defining $g^{(d)}$ and such that ${D_{\phi}}(g^{(d)},f)=0$ exist, then $$\frac{\dr^d}{\dr y_1...\dr y_d}\tilde C_f=\frac{\dr^d}{\dr y_1...\dr y_d}\tilde C_g$$
\end{proposition}
\subsection{With the elliptical copula}\label{CopulaTest}

Let $f$ be an unknown density defined on $\R^d$. The objective of the present section is to determine whether the copula of $f$ is elliptical. We thus define an instrumental elliptical density $g$ with the same mean and variance as  $f$, and we follow the procedure of section \ref{UseSample}. As explained in section \ref{ExmapleCopula}, we infer from proposition \ref{PropCop} that the copula of $f$ equals the copula of $g$ when ${D_{\phi}}(g^{(d)},f)=0$, i.e. when $a_d$ is the last vector generated from the algorithm and when $(a_i)$ is the canonical basis of $\R^d$. Thus, in order to verify this assertion, corollary \ref{KernelLOIDUCRITERE2} page \pageref{KernelLOIDUCRITERE2} provides us with a $\alpha$-level confidence ellipsoid around this vector, namely 
$$\cE_d=\{b\in\R^d;\ \sqrt n(Var_{\PP}(M(b,b)))^{-1/2}\Pn_nM(b,b)\leq q_{\alpha}^{\cN(0,1)} \}$$
where $q_{\alpha}^{\cN(0,1)}$ is the quantile of a $\alpha$-level reduced centered normal distribution, where $\Pn_n$ is the empirical measure araising from a realization of the sequences $(X_1,\ldots,X_n)$ and $(Y_1,\ldots,Y_n)$ - see appendix \ref{truncSample} - and where $M$ is a known function of $d$ , $f_n$ and $g^{(d-1)}_n$ - see section \ref{HypoF}.\\
Consequently, keeping the notations introduced in section \ref{ExmapleCopula}, we can perform a statistical test of the null hypothesis
$$ \text{ $(H_0)$ : $\frac{\dr^d}{\dr x_1...\dr x_d}C_f=\frac{\dr^d}{\dr x_1...\dr x_d} C_g$ versus 
$(H_1)$ : $\frac{\dr^d}{\dr x_1...\dr x_d} C_f\not=\frac{\dr^d}{\dr x_1...\dr x_d} C_g$}$$
Since, under $(H_0)$, we have ${D_{\phi}}(g^{(d)},f)=0$, then the following theorem provides us with 
a confidence region for this test.
\begin{theoreme}\label{TestCop}
The set $\cE_d$ is a confidence region for the test of the null hypothesis $(H_0)$ versus the alternative $(H_1)$.\end{theoreme}
\begin{remarque}
1/ If ${D_{\phi}}(g^{(k)},f)=0$, for $k<d$, then we reiterate the algorithm until $g^{(d)}$ is created in order to obtain a relationship for the copula of $f$.\\
2/ If the $a_i$ do not constitute the canonical basis, then keeping the notations introduced in section \ref{ExmapleCopula}, our algorithm meets the test :
$$ \text{ $(H_0)$ : $\frac{\dr^d}{\dr y_1...\dr y_d}\tilde C_f=\frac{\dr^d}{\dr y_1...\dr y_d}\tilde C_g$ versus 
$(H_1)$ : $\frac{\dr^d}{\dr y_1...\dr y_d}\tilde C_f\not=\frac{\dr^d}{\dr y_1...\dr y_d}\tilde C_g$}$$
Thus, our method enables us to tell wether the copula of $f$ equals the copula of $g$ in the $(a_1,\ldots,a_d)$ basis.
\end{remarque}
\subsection{With the independent copulas}\label{CopulaIndepTest}
Let $f$ be a density on $\R^d$ and let $X$ be a random vector with $f$ as density. The objective of this section is to determine whether $f$ is the product of its margins, i.e. whether the copula of $f$ is the independent copula. Let thus $g$ be an instrumental product of univariate Gaussian density - with $diag(Var(X_1),...,Var(X_d))$ as covariance matrix and with the same mean as $f$ - as explained at section \ref{CopulaTest}, let us follow the procedure described at section \ref{UseSample}, i.e. proposition \ref{PropCop} infers that the copula of $f$ is the independent copula when ${D_{\phi}}(g^{(d)},f)=0$. Thus, we perform a statistical test of the null hypothesis :
$$ \text{ $(H_0)$ : $f=\Pi_{i=1}^df_i$ versus the alternative $(H_1)$ : $f\not=\Pi_{i=1}^df_i$} $$ 
Since, under $(H_0)$, we have ${D_{\phi}}(g^{(d)},f)=0$, then the following theorem provides us with a confidence region for our test.
\begin{theoreme}\label{TestIndepCop}
Keeping the notations of section \ref{CopulaTest}, the set $\cE_d$ is a confidence region for the test of the null hypothesis $(H_0)$ versus the alternative $(H_1)$.
\end{theoreme}
\begin{remarque}
1/ As explained in section \ref{CopulaTest}, if ${D_{\phi}}(g^{(k)},f)=0$, for $k<d$, we reiterate the algorithm until $g^{(d)}$ is created in order to derive a relationship for the copula of $f$.\\
2/ If the $a_i$ do not constitute the canonical basis, then keeping the notations introduced in section \ref{ExmapleCopula}, our algorithm meets the test :
$$ \text{ $(H_0)$ : $f=\Pi_{i=1}^df_{a_i}$ versus the alternative $(H_1)$ : $f\not=\Pi_{i=1}^df_{a_i}$} $$ 
Thus, our method enables us to determine if the the copula of $f$ is the independent copula in the $(a_1,\ldots,a_d)$ basis.
\end{remarque}
\subsection{Study of the subsequence $(g^{(k')})$ defined by ${D_{\phi}}(g^{(k')},f)=0$ for any $k'$}\label{SousSuiteSection}
Let $\cQ$ be the set of non-negative integers defined by $\cQ=\{k_i';\ k_1'=1,\ k_q'=d,\ k_i'<k_{i+1}'\}$, where $q$ - such that $q\leq d$ - is its cardinal. In the present section, our goal is to study the subsequence $(g^{(k')})$ of the sequence $(g^{(k)})_{k=1..d}$ defined by ${D_{\phi}}(g^{(k')},f)=0$ for any $k'$ belonging to $\cQ$.\\
First, we have :\\
${D_{\phi}}(g^{(d)},f)=0\ \Leftrightarrow\ g^{(d)}=f$, through property \ref{Phimini}

$\Leftrightarrow\ \frac{g(x)}{\Pi_{k=1}^d g_{a_k}(\transp {a_k}x)}=\frac{f(x)}{\Pi_{k=1}^df_{a_k}(\transp {a_k}x)}$, as explained in section \ref{CopulaTest},

$\Leftrightarrow\ \frac{\tilde g(y)}{\Pi_{k=1}^d \tilde g_{k}(y_k)}=\frac{\tilde f(y)}{\Pi_{k=1}^d\tilde f_{k}(y_k)}$, which amounts to the previous relationship written in the $A=(a_1,\ldots,a_d)$ 

basis with the notations introduced in section \ref{CopulaTest}.\\
Moreover, defining $\tilde k_i'$ as the previous integer $k_i'$, in the space $\{1,\ldots,d\}$, with $i>1$, and as explained in section \ref{modelSection}, the relationship ${D_{\phi}}(g^{(k')},f)=0$ implies that

$\tilde f(y_i,\ldots,y_{\tilde k_{i+1}'}/y_1,\ldots,y_{\tilde k_{i}'},y_{k_{i+1}'},\ldots,y_{d})=\tilde f_{i,i+1}(y_i,\ldots,y_{\tilde k_{i+1}'})$\\
where $\tilde f_{i,i+1}$ is the density of vector $(\transp{a_i}X,\ldots,\transp{a_{\tilde k_{i+1}'}}X)$ in the $A=(a_1,\ldots,a_d)$ basis.\\ Consequently,
$\tilde f(y)=\tilde f_{1,2}(y_1,\ldots,y_{\tilde k_2'}).\tilde f_{2,3}(y_{k_2'},\ldots,y_{\tilde k_3'})\ldots\tilde f_{q-1,d}(y_{k_{q-1}'},\ldots,y_{\tilde k_d'})$.\\
Hence, we can infer that
\begin{equation}\label{SuperTruc}
\frac{\tilde f(y)}{\Pi_{k=1}^d\tilde f_{k}(y_k)}=
\frac{\tilde f_{1,2}(y_1,\ldots,y_{\tilde k_2'})}
     {\Pi_{k=1}^{\tilde k_2'}\tilde f_{k}(y_k)}.
\frac{\tilde f_{2,3}(y_{k_2'},\ldots,y_{\tilde k_3'})}
     {\Pi_{k={k_2'}}^{\tilde k_3'}\tilde f_{k}(y_k)}
\ldots
     \frac{\tilde f_{q-1,d}(y_{k_{q-1}'},\ldots,y_{\tilde k_d'})}
     {\Pi_{k={\tilde k_{q-1}'}}^d\tilde f_{k}(y_k)}
\end{equation}
The following theorem explicitely describes the form of the $f$ copula in the $A=(a_1,\ldots,a_d)$ basis :
\begin{theoreme}\label{Theo-Real-Case}
Defining $\tilde C_{f_{i,j}}$ as the copula of $\tilde f_{i,j}$ and keeping the notations introduced in sections \ref{ExmapleCopula} and \ref{SousSuiteSection}, it holds
$$
\frac{\dr^d}{\dr y_1...\dr y_d}\tilde C_f=
\frac{\dr^{\tilde k_2'}}{\dr y_1\ldots\dr y_{\tilde k_2'}}\tilde C_{f_{1,2}}.
\frac{\dr^{\tilde k_3'-k_2'+1}}{\dr y_{k_2'}\ldots\dr y_{\tilde k_3'}}\tilde C_{f_{2,3}}
\ldots
\frac{\dr^{d-k_{q-1}'+1}}{\dr y_{k_{q-1}'}\ldots\dr y_{d}}\tilde C_{f_{q-1,d}}$$
\end{theoreme}
\begin{remarque}\label{RemCopule=1}
If there exists $i$ such that $i<d$ and $k_i'=\tilde k_{i+1}'$, then the notation $\tilde f_{i,i+1}(y_{k_{i}'},\ldots,y_{\tilde k_{i+1}'})$ means $\tilde f_{k_{i}'}(y_{k_{i}'})$. Thus, if, for any $k$, we have ${D_{\phi}}(g^{(k)},f)=0$, then, for any $i<d$, we have $k_i'=\tilde k_{i+1}'$, i.e. we have $\tilde f=\Pi_{k=1}^d\tilde f_k(y_k)$ - where $\tilde f_k$ is the $k^{th}$ marginal density of $\tilde f$.
\end{remarque}
At present, using relationship \ref{SuperTruc} and remark \ref{RemCopule=1}, the following corollary gives us the copula of $f$ as equals to 1 in the $\{a_1,\ldots,a_d\}$ basis when, for any $k$, ${D_{\phi}}(g^{(k')},f)=0$ :
\begin{corollaire}\label{Coro-Real-Case}
In the case where, for any $k$, ${D_{\phi}}(g^{(k)},f)=0$, it holds:
$$
\frac{\dr^d}{\dr y_1...\dr y_d}\tilde C_f=1
$$
\end{corollaire}
\section{ Simulations}\label{Simul400}

Let us examine two simulations and an application to real datasets. The first simulation studies the elliptical copula and the second studies the independent copula.
In each simulation, our program will aim at creating a sequence of densities $(g^{(j)})$, $j=1,..,d$ such that $g^{(0)}=g,$ $g^{(j)}=g^{(j-1)}f_{a_j}/[g^{(j-1)}]_{a_j}$ and ${D_{\phi}}(g^{(d)},f)=0,$ where ${D_{\phi}}$ is a divergence - see appendix \ref{FiDiv} for its definition - and\\ $a_j=arg\inf_b {D_{\phi}}(g^{(j-1)}f_b/g^{(j-1)}_b,f),$ for all $j=1,...,d$. We will therefore perform the tests introduced in theorems \ref{TestCop} and \ref{TestIndepCop}.
\begin{simulation}\label{Sim1}\rm
$\\$We are in dimension $2$(=d), and we use the $\chi^2$ divergence to perform our optimisations. Let us consider a sample of $50$(=n) values of a random variable $X$ with a density law $f$ defined by :
$$f(x)=c_\rho(F_{Gumbel}(x_1),F_{Exponential}(x_2)).Gumbel(x_1).Exponential(x_2)$$
where :\\
$\bullet$ $c$ is the Gaussian copula with correlation coefficient $\rho=0.5$,\\
$\bullet$ the  Gumbel distribution parameters are $-1$ and $1$ and the exponential density parameter is $2$. \\
Let us generate then a Gaussian random variable $Y$ with a density - that we will name $g$ - presenting the same mean and variance as $f$.\\
We theoretically obtain $k=2$ and $(a_1,a_2)=((1,0),(0,1))$.\\
To get  this result, we perform the following test: 
$$(H_0):\ (a_1,a_2)=((1,0),(0,1))\text{ versus }(H_1):\ (a_1,a_2)\not=((1,0),(0,1))$$
Then, theorem \ref{TestCop} enables us to verify $(H_0)$ by the following 0.9(=$\alpha$) level confidence ellipsoid
$${\mathcal E}_{2}=\{b \in \R^2;\ (Var_{\bf P}(M(b,b)))^{(-1/2)}\Pn_nM(b,b)\leq q^{\cN(0,1)}_{\alpha}/\sqrt n\simeq 0,2533/7.0710=0.03582\}$$
And, we obtain
\begin{table}[h]
\caption{Simulation 1 : Numerical results of the optimisation}
\label{Sim1-opti}
\begin{tabular}{cll}
\hline
Our Algorithm &&\\
\hline
\multirow{3}{8cm}{Projection Study 0 :}
 & minimum :  0.445199        \\  \cline{2-3}
 & at point : (1.0171,0.0055) \\  \cline{2-3}
 & P-Value :  0.94579         \\
\hline
\multirow{1}{8cm}{Test :}
 & $H_1$ : $a_1\not\in {\mathcal E}_{1}$ : True  \\
\hline
\multirow{3}{8cm}{Projection Study 1 :}
 & minimum :  0.009628        \\  \cline{2-3}
 & at point : (0.0048,0.9197) \\  \cline{2-3}
 & P-Value :  0.99801         \\
\hline
\multirow{1}{8cm}{Test :}
 & $H_0$ : $a_2\in {\mathcal E}_{2}$ : True   \\
\hline
\multirow{1}{8cm}{$\chi^2$(Kernel Estimation of $g^{(2)}$, $g^{(2)}$)}
 & 3.57809  \\
\hline
\end{tabular}

Therefore, we can conclude that $H_0$ is verified.
\end{table}
\begin{figure}[h!]
  \center
  \caption{Graph of the estimate of $(x_1,x_2)\mapsto c_\rho(F_{Gumbel}(x_1),F_{Exponential}(x_2))$.}
  \includegraphics[scale=0.45]{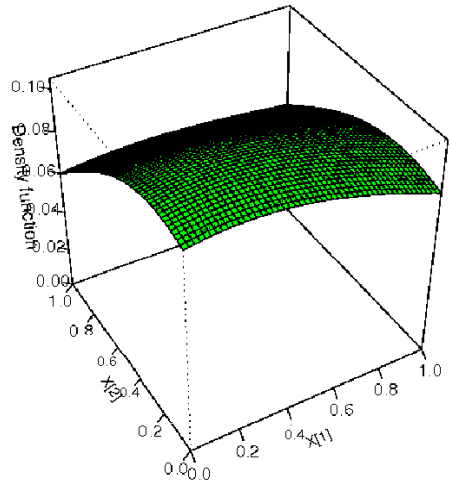}
\end{figure}
\end{simulation}
\newpage
\begin{simulation}\label{Sim3}\rm
$\\$We are in dimension $2$(=d), and we use the $\chi^2$ divergence to perform our optimisations.\\
Let us consider a sample of $50$(=n) values of a random variable $X$ with a density law $f$ defined by 

$f(x)=Gumbel(x_1).Exponential(x_2)$,\\
where the  Gumbel distribution parameters are $-1$ and $1$ and the exponential density parameter is $2$.\\
Let $g$ be an instrumental product of univariate Gaussian densities - with $diag(Var(X_1),...,Var(X_d))$ as covariance matrix and with the same mean as $f$.\\
We theoretically obtain $k=2$ and $(a_1,a_2)=((1,0),(0,1))$. To get  this result, we perform the following test:

$(H_0):\ (a_1,a_2)=((1,0),(0,1))\text{ versus }(H_1):\ (a_1,a_2)\not=((1,0),(0,1))$.\\
Then, theorem \ref{TestIndepCop} enables us to verify $(H_0)$ by the following 0.9(=$\alpha$) level confidence ellipsoid 

${\mathcal E}_{2}=\{b \in \R^2;\  (Var_{\bf P}(M(b,b)))^{(-1/2)}\Pn_nM(b,b)\leq q^{\cN(0,1)}_{\alpha}/\sqrt n\simeq 0.03582203\}$.\\
And, we obtain

\begin{table}[h]
\caption{Simulation 2 : Numerical results of the optimisation}
\label{Sim2-opti}
\begin{tabular}{cll}
\hline
Our Algorithm &&\\
\hline
\multirow{3}{8cm}{Projection Study 0 :}
 & minimum :  0.057833        \\  \cline{2-3}
 & at point : (0.9890,0.1009) \\  \cline{2-3}
 & P-Value :  0.955651        \\
\hline
\multirow{1}{8cm}{Test :}
 & $H_1$ : $a_1\not\in {\mathcal E}_{1}$ : True  \\
\hline
\multirow{3}{8cm}{Projection Study 1 :}
 & minimum :  0.02611          \\  \cline{2-3}
 & at point : (-0.1105,0.9290) \\  \cline{2-3}
 & P-Value :  0.921101         \\
\hline
\multirow{1}{8cm}{Test :}
 & $H_0$ : $a_2\in {\mathcal E}_{2}$ : True   \\
\hline
\multirow{1}{8cm}{$\chi^2$(Kernel Estimation of $g^{(2)}$, $g^{(2)}$)}
 & 1.25945 \\
\hline
\end{tabular}

Therefore, we can conclude that $f=\Pi_{i=1}^df_i.$
\end{table}
\begin{figure}[h!]
 \center
 \caption{Graph of the independent copula estimate.}
 \includegraphics[scale=0.5]{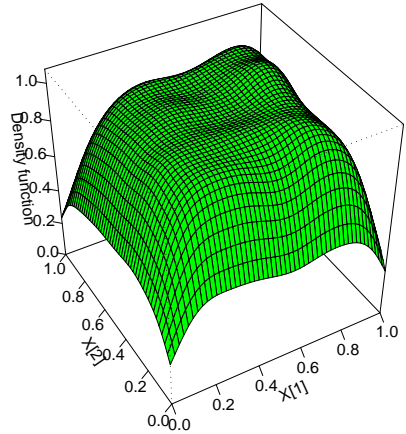}
\end{figure}
\end{simulation}
\newpage
\subsubsection{\bf Application to real datasets}
Let us for instance study the moves in the stock prices of Renault and Peugeot from January 4, 2010 to July 25, 2010. We thus gather 140(=n) data from these stock prices - see data below.\\
Let us also consider $X_1$ (resp. $X_2$) the random variable defining the stock price of Renault (resp. Peugeot). We will assume - as it is commonly done in mathematical finance - that the stock market abides by the classical hypotheses of the Black-Scholes model - see \citet{BS73}.\\
Consequently, $X_1$ and $X_2$ each present a log-normal distribution as probability distribution.\\
Let $f$ be the density of vector $(ln(X_1),ln(X_2))$, let us now apply our algorithm to $f$ with the Kullback-Leibler divergence as $\phi$-divergence. Let us generate then a Gaussian random variable $Y$ with a density - that we will name $g$ - presenting the same mean and variance as $f$.\\
We first assume that there exists a vector $a$ such that ${D_{\phi}}(g\frac{f_a}{g_a},f)=0$.\\
In order to verify this hypothesis, our reasoning will be the same as in Simulation \ref{Sim1}. Indeed, we assume that this vector is a co-factor of $f$. Consequently, corollary \ref{KernelLOIDUCRITERE2coro} enables us to estimate $a$ by the following 0.9(=$\alpha$) level confidence ellipsoid 

${\mathcal E}_{1}=\{b \in \R^2;\  (Var_{\bf P}(M(b,b)))^{(-1/2)}\Pn_nM(b,b)\leq q^{\cN(0,1)}_{\alpha}/\sqrt n \simeq0,2533/\sqrt{140}=0.02140776\}$.\\ And, we obtain
\begin{table}[h]
\caption{Numerical results : First projection}
\label{Real-opti}
\begin{tabular}{cll}
\hline
Our Algorithm &&                                      \\
\hline
\multirow{3}{8cm}{Projection Study 0 :}
 & minimum   :  0.02087685                            \\  \cline{2-3}
 & at point  : $a_1$=(19.1,-12.3)                     \\  \cline{2-3}
 & P-Value   :  0.748765                              \\
\hline
\multirow{1}{8cm}{Test :}
 & $H_0$     : $a_1\in {\mathcal E}_{1}$ : True       \\
\hline
\multirow{1}{8cm}{K(Kernel Estimation of $g^{(1)}$, $g^{(1)}$)}
& 4.3428735                                           \\
\hline
\end{tabular}
\end{table}
$\\$Therefore, our first  hypothesis is confirmed.\\
However, our goal is to study the copula of $(ln(X_1),ln(X_2))$. Then, as explained in section \ref{SousSuiteSection}, we formulate another hypothesis assuming that there exists a vector $a$ such that ${D_{\phi}}(g^{(1)}\frac{f_a}{g_a^{(1)}},f)=0$.\\
In order to verify this hypothesis, we will use the same reasoning as above. Indeed, we assume that this vector is a co-factor of $f$. Consequently, corollary \ref{KernelLOIDUCRITERE2coro} enables us to estimate $a$ by the following 0.9(=$\alpha$) level confidence ellipsoid 

${\mathcal E}_{2}=\{b \in \R^2;\  (Var_{\bf P}(M(b,b)))^{(-1/2)}\Pn_nM(b,b)\leq q^{\cN(0,1)}_{\alpha}/\sqrt n \simeq0,2533/\sqrt{140}=0.02140776\}$.\\ And, we obtain
\begin{table}[h]
\caption{Numerical results : Second projection}
\label{Real-opti2}
\begin{tabular}{cll}
\hline
Our Algorithm &&                                      \\
\hline
\multirow{3}{8cm}{Projection Study 1 :}
 & minimum   :  0.0198753                             \\  \cline{2-3}
 & at point  : $a_2$=(8.1,3.9)                        \\  \cline{2-3}
 & P-Value   :  0.8743401                             \\
\hline
\multirow{1}{8cm}{Test :}
 & $H_0$     : $a_2\in {\mathcal E}_{2}$ : True       \\
\hline
\multirow{1}{8cm}{K(Kernel Estimation of $g^{(2)}$, $g^{(2)}$)}
& 4.38475324                                          \\
\hline
\end{tabular}
\end{table}
$\\$Therefore, our second  hypothesis is confirmed.

$\text{In conclusion, as explained in corollary \ref{Coro-Real-Case}, the copula of $f$ is equal to $1$ in the $\{a_1,a_2\}$ basis.}$
\begin{figure}[h!]
 \center
 \caption{Graph of the copula of $(ln(X_1),ln(X_2))$ in the canonical basis.}
 \includegraphics[scale=0.45]{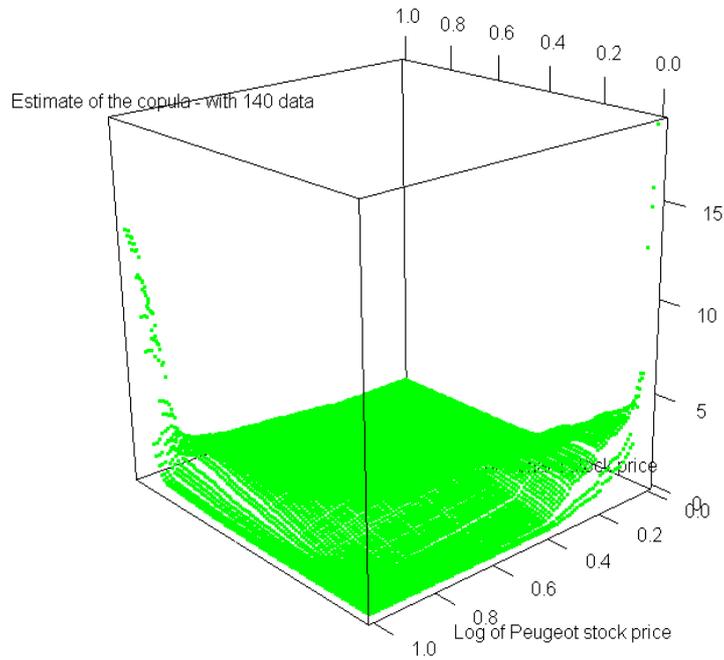}
 \label{real-case}
\end{figure}
\begin{figure}[h!]
 \center
 \caption{Graph of the copula of $(ln(X_1),ln(X_2))$ in the $\{a_1,a_2\}$ basis.}
 \includegraphics[scale=0.45]{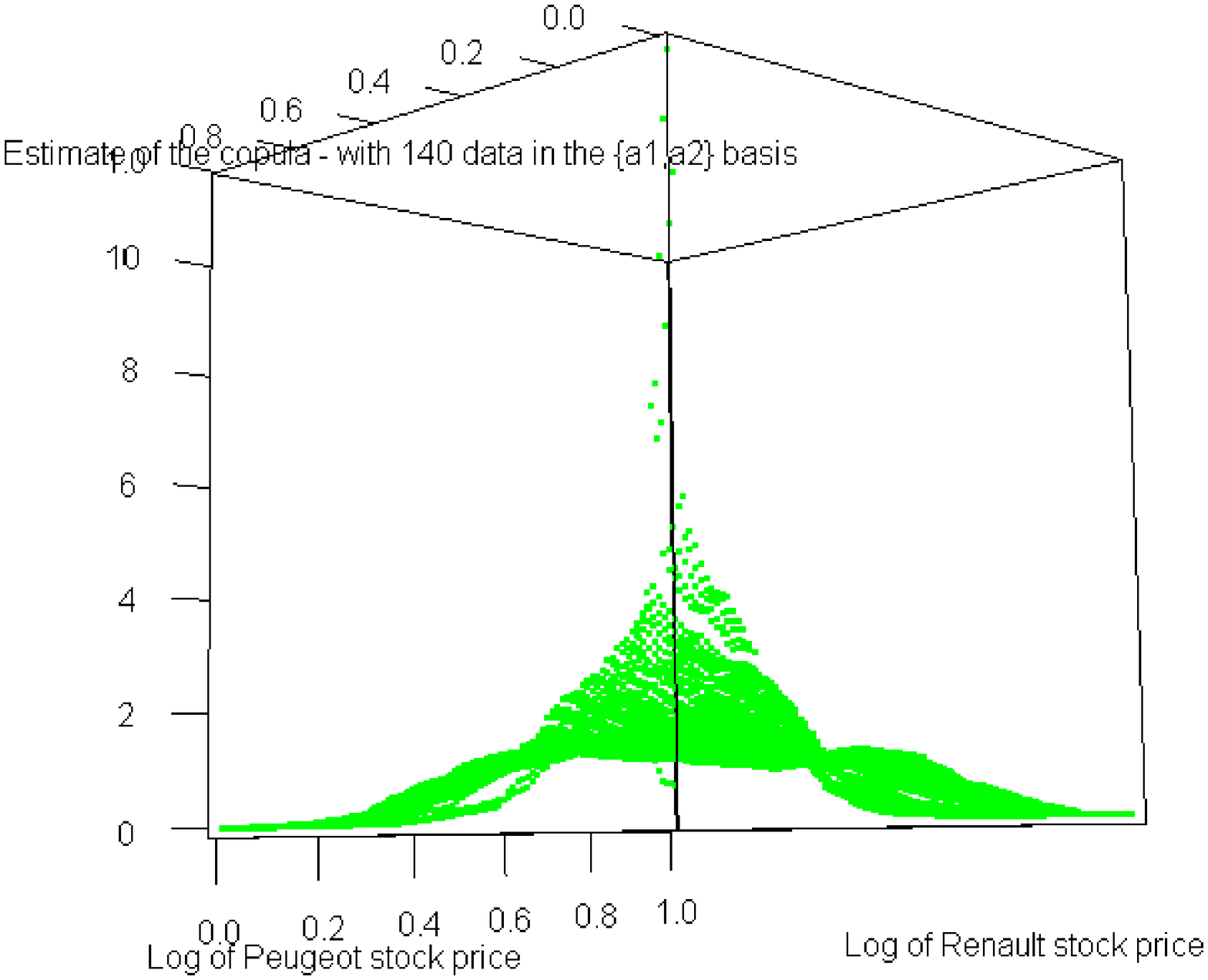}
 \label{real-case-ab-1}
\end{figure}
\begin{figure}[h!]
 \center
 \caption{Graph of the copula of $(ln(X_1),ln(X_2))$ in the $\{a_1,a_2\}$ basis - other view.}
 \includegraphics[scale=0.45]{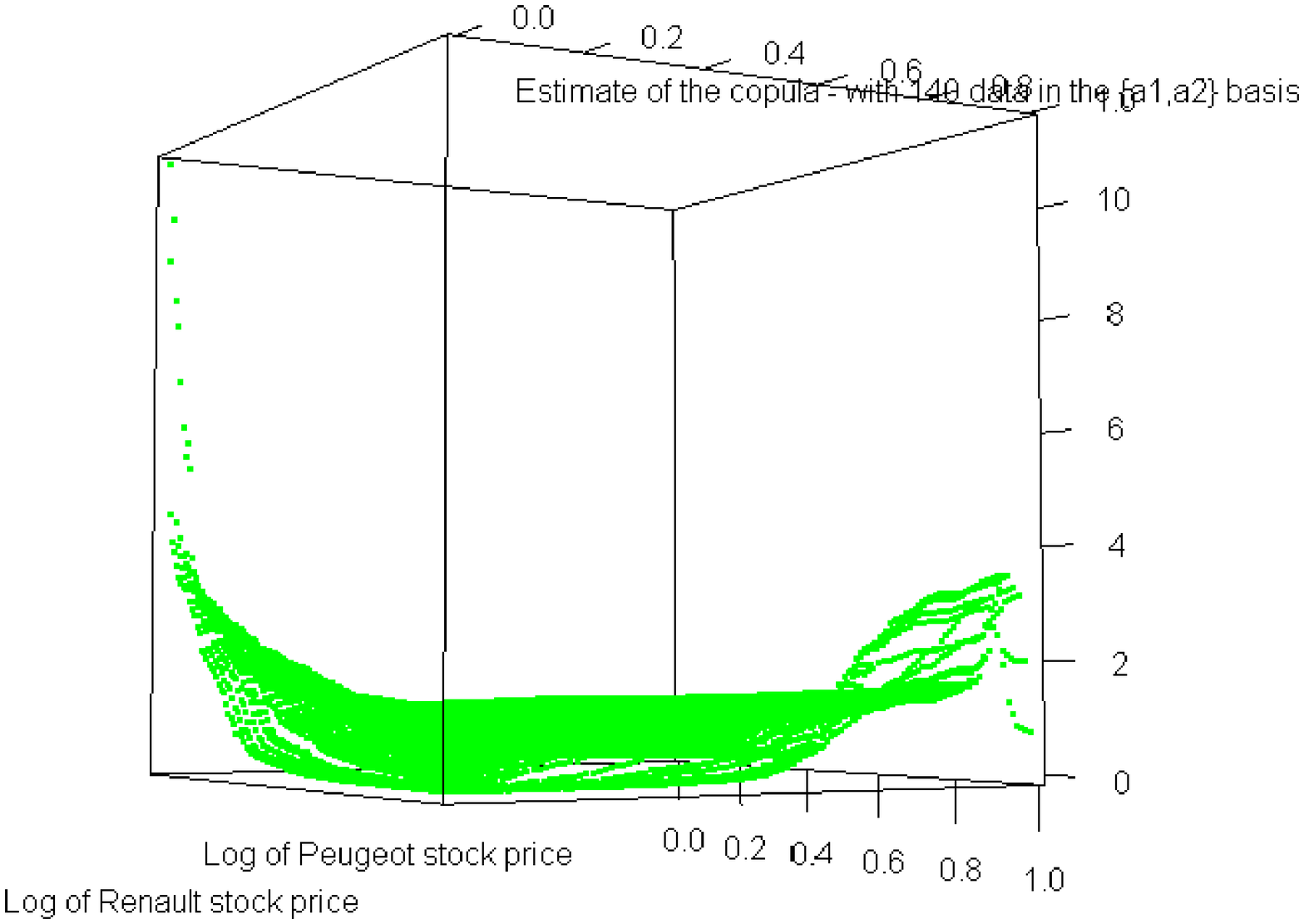}
 \label{real-case-ab-2}
\end{figure}
\begin{table}[htb!]
\caption{Stock prices of Renault and Peugeot}
\label{Real-data2}
\begin{tabular}{|c|c|c||c|c|c||c|c|c|}
\hline
Date&Renault&Peugeot&Date&Renault&Peugeot&Date&Renault&Peugeot\\
\hline
23/07/10&34.9&24.2&22/07/10&34.26&24.01&21/07/10&33.15&23.3\\
20/07/10&32.69&22.78&19/07/10&33.24&23.36&16/07/10&33.92&23.77\\
15/07/10&34.44&23.71&14/07/10&35.08&24.36&13/07/10&35.28&24.37\\
12/07/10&33.84&23.16&09/07/10&33.46&23.13&08/07/10&33.08&22.65\\
07/07/10&32.15&22.19&06/07/10&31.12&21.56&05/07/10&30.02&20.81\\
02/07/10&30.17&20.85&01/07/10&29.56&20.05&30/06/10&30.78&21.07\\
29/06/10&30.55&20.97&28/06/10&32.34&22.3&25/06/10&31.35&21.68\\
24/06/10&32.29&22.25&23/06/10&33.58&22.47&22/06/10&33.84&22.77\\
21/06/10&34.06&23.25&18/06/10&32.89&22.7&17/06/10&32.08&22.31\\
16/06/10&31.87&21.92&15/06/10&32.03&22.12&14/06/10&31.45&22.2\\
11/06/10&30.62&21.42&10/06/10&30.42&20.93&09/06/10&29.27&20.34\\
08/06/10&28.48&19.73&07/06/10&28.92&20.15&04/06/10&29.19&20.27\\
03/06/10&30.35&20.46&02/06/10&29.33&19.53&01/06/10&28.87&19.45\\
31/05/10&29.39&19.54&28/05/10&29.16&19.55&27/05/10&29.18&19.81\\
26/05/10&27.5&18.5&25/05/10&26.76&18.08&24/05/10&28.75&18.81\\
21/05/10&28.78&18.82&20/05/10&28.53&18.84&19/05/10&29.49&19.25\\
18/05/10&30.95&19.76&17/05/10&30.92&19.35&14/05/10&31.35&19.34\\
13/05/10&33.65&20.76&12/05/10&33.63&20.52&11/05/10&33.38&20.34\\
10/05/10&33.28&20.3&07/05/10&31&19.24&06/05/10&32.4&20.22\\
05/05/10&32.95&20.45&04/05/10&33.3&21.03&03/05/10&35.58&22.63\\
30/04/10&35.41&22.45&29/04/10&35.53&22.36&28/04/10&34.75&22.33\\
\hline
\end{tabular}
\end{table}
\begin{table}[htb!]
\caption{Stock prices of Renault and Peugeot (continued)}
\label{Real-data3}
\begin{tabular}{|c|c|c||c|c|c||c|c|c|}
\hline
Date&Renault&Peugeot&Date&Renault&Peugeot&Date&Renault&Peugeot\\
\hline
27/04/10&36.2&22.9&26/04/10&37.65&23.73&23/04/10&36.72&23.5\\
22/04/10&34.36&22.72&21/04/10&35.01&22.86&20/04/10&35.62&22.88\\
19/04/10&34.08&21.77&16/04/10&34.46&21.71&15/04/10&35.16&22.22\\
14/04/10&35.1&22.22&13/04/10&35.28&22.45&12/04/10&35.17&21.85\\
09/04/10&35.76&21.9&08/04/10&35.67&21.67&07/04/10&36.5&21.89\\
06/04/10&36.87&22&01/04/10&35.5&21.97&31/03/10&34.7&21.8\\
30/03/10&34.8&22.24&29/03/10&35.7&22.73&26/03/10&35.54&22.58\\
25/03/10&35.53&22.73&24/03/10&33.8&21.82&23/03/10&34.1&21.58\\
22/03/10&33.73&21.64&19/03/10&34.12&21.68&18/03/10&34.44&21.75\\
17/03/10&34.68&21.98&16/03/10&34.33&21.88&15/03/10&33.57&21.53\\
12/03/10&33.9&21.86&11/03/10&33.27&21.58&10/03/10&33.12&21.47\\
09/03/10&32.69&21.54&08/03/10&32.99&21.66&05/03/10&32.89&21.85\\
04/03/10&31.64&21.26&03/03/10&31.65&20.7&02/03/10&31.05&20.2\\
01/03/10&30.26&19.54&26/02/10&30.2&19.39&25/02/10&29.42&18.98\\
24/02/10&30.9&19.49&23/02/10&30.54&19.74&22/02/10&31.89&20.06\\
19/02/10&32.29&20.67&18/02/10&32.26&20.41&17/02/10&31.69&20.31\\
16/02/10&31.08&19.8&15/02/10&30.25&19.66&12/02/10&29.56&19.57\\
11/02/10&31&20.4&10/02/10&32.78&21.21&09/02/10&33.31&22.31\\
08/02/10&32.63&21.95&05/02/10&32.15&22.33&04/02/10&33.72&22.86\\
03/02/10&35.32&23.93&02/02/10&35.29&23.8&01/02/10&35.31&24.05\\
29/01/10&34.26&23.64&28/01/10&33.94&23.31&27/01/10&33.85&23.88\\
26/01/10&34.97&24.86&25/01/10&35.06&24.35&22/01/10&35.7&24.95\\
21/01/10&36.1&25&20/01/10&36.92&25.35&19/01/10&38.4&25.81\\
18/01/10&39.28&25.95&15/01/10&38.6&25.7&14/01/10&39.56&26.67\\
13/01/10&39.49&26.13&12/01/10&38.36&25.98&11/01/10&39.21&26.65\\
08/01/10&39.38&26.5&07/01/10&39.69&26.7&06/01/10&39.25&26.32\\
05/01/10&38.31&24.74&04/01/10&38.2&24.52&       &     &     \\
\hline
\end{tabular}
\end{table}

\vspace{100mm}

\newpage
$\text{\vspace{100mm}}$\\

\subsection*{\bf Critics of the simulations}
In the case where $f$ is unknown, we will never be sure to have reached the minimum of the $\phi$-divergence:  we have indeed used the simulated annealing method to solve our optimisation problem, and therefore it is only when the number of random jumps tends in theory towards infinity that the probability to get the minimum tends to 1.
We also note that no theory on the optimal number of jumps to implement does exist, as this number depends on the specificities of each particular problem.\\
Moreover, we choose the $50^{-\frac{4}{4+d}}$ for the AMISE of the two simulations. This choice leads us to simulate 50 random variables - see \citet{MR1191168} page 151 -, none of which have been discarded to obtain the truncated sample.\\
This has also been the case in our application to real datasets.\\
Finally, the shape of the copula in the case of real datasets in the $\{a_1,\ a_2\}$ basis is also noteworthy. \\
Figure \ref{real-case-ab-1} shows  that the curve reaches a quite wide plateau around $1$, whereas Figure \ref{real-case-ab-2} shows that this plateau prevails on almost the entire $[0,1]^2$ set. We can therefore conclude that the theoritical analysis is indeed confirmed by the above simulation.

\subsection*{\bf  Conclusion}
Projection Pursuit is useful in evidencing characteristic structures as well as one-dimensional projections and their associated distribution in multivariate data. This article clearly evidences the efficiency of the $\varphi$-projection pursuit methodology for goodness-of-fit tests for copulas. Indeed, the robustness as well as the convergence results we achieved, convincingly fulfilled our expectations regarding the methodology used.
\appendix
\section{On the different families of copula}\label{CopFam}

There exists many copula families. Let us here present the most important amongst them.

\subsection{Archimedean copulas}
These copulas present  a simple form with properties such as associativity and have a variety of dependent structures. They can generally be defined under the following form
$$ H(u_1,u_2,\dots,u_n) = \Psi^{-1}\left(\sum_{i=1}^n\Psi(F_i(u_i)) \right)$$
where $(u_1,u_2,\dots,u_n)\in[0,1]^n$ and where $\Psi$ is known as a "generator function". This $\Psi$ function must be at least $d-2$ times continuously differentiable, must have a decreasing and convex $d-2$ derivative, and must be such that $\Psi(1)=0$.

Let us now present several examples :

{\it 1/ Clayton copula:}\\
The Clayton copula is an asymmetric archimedean copula, exhibiting greater dependency in the negative tail than in the positive tail.
Let us define $X$ (resp. $Y$) as the random vector having $F$ (resp $G$) as cumulative distribution function (CDF). Assuming that the vector $(X,Y)$ has a Clayton copula, then this copula is given by:
$$H(x,y) = (F(x)^\theta+G(y)^\theta-1)^{1/\theta}$$
And its generator is:
$$\Psi(x) = x^{\theta} -1$$
For $\theta= 0$ in the Clayton copula, the random variables are statistically independent. 
The generator function approach can be extended to create multivariate copulas, simply by including more additive terms.

{\it 2/ Gumbel copula:}\\
The Gumbel copula (a.k.a. Gumbel-Hougard copula) is an asymmetric archimedean copula, exhibiting greater dependency in the positive tail than in the negative tail. This copula is given by:
$$\Psi(x)= (-\ln(x))^\alpha$$
{\it 3/ Frank copula:}\\
The Frank copula is a symmetric archimedean copula given by:
$$\Psi(x)= \ln\left( \frac{e^{\alpha x} -1}{e^{\alpha} -1}\right)$$
\subsection{Periodic copula}
In 2005, \citet{PeriodicCopula} introduced a way of constructing copulas based on periodic functions.
Defining $h$ (resp. $H$) as a 1-periodic non-negative function that integrates to $1$ over $[\ 0,\ 1]$ (resp. as a double primitive of $h$), then both
$$H(u+v)-H(u)-H(v)\text{ and }-H(u-v)+H(u)+H(-v)$$
are copula functions, the second one not being necessarily exchangeable.
\section{ $\phi$-Divergence}\label{FiDiv}

Let us call $h_a$ the density of  $\transp a Z$ if  $h$ is the density of $Z$. 
Let $\varphi$ be a strictly convex function defined by $\varphi: \overline{\R^+}\to\overline{\R^+},$ and such that $\varphi(1)=0$.
\begin{definition}
We define a $\phi-$divergence of $P$ from $Q$, where $P$ and $Q$ are two probability distributions over a space $\Omega$ such that $Q$ is absolutely continuous with respect to $P$, by
\begin{equation}\label{def:div}
{D_{\phi}}(Q,P)=\int\varphi(\frac{dQ}{dP})dP.
\end{equation}
The above expression (\ref{def:div}) is also valid if $P$ and $Q$ are both dominated by the same probability.
\end{definition}
The most used distances (Kullback, Hellinger or $ \chi^2$) belong to the Cressie-Read family (see \citet{CressieRead}, \citet{Csiszar67} and the books of \citet{MR926905}, \citet{MR2183173} and \citet{MR1075502}). They are defined by a specific $\varphi$. Indeed,\\
- with the Kullback-Leibler divergence, we associate $\varphi(x)=xln(x)-x+1$\\
- with the Hellinger distance, we associate $\varphi(x)=2(\sqrt x-1)^2$\\ 
- with the $\chi^2$ distance, we associate $\varphi(x)=\frac{1}{2}(x-1)^2$\\
- more generally, with power divergences, we associate $\varphi(x)=\frac{x^\gamma-\gamma x+\gamma-1}{\gamma(\gamma-1)}$, where $\gamma\in\R\setminus {(0,1)}$\\
- and, finally, with the $L^1$ norm, which is also a divergence, we associate $\varphi(x)=|x-1|.$\\
Let us now present some well-known properties of divergences.
\begin{PROPIE}\label{Phimini}
We have ${D_{\phi}}(P,Q)=0\Leftrightarrow P=Q.$
\end{PROPIE}
\begin{PROPIE}\label{K-SCI}
The divergence function $Q\mapsto {D_{\phi}}(Q,P)$ is\\
$\bullet$  convex,\\
$\bullet$  lower semi-continuous, for the topology that makes all the applications of the form $Q\mapsto\int fdQ$ continuous  where $f$ is bounded and continuous, and\\
$\bullet$ lower semi-continuous for the topology of the uniform convergence.
\end{PROPIE}
Finally, we will also use the following property derived from the first part of corollary  (1.29) page 19 of \citet{MR926905},
\begin{PROPIE}\label{ExitenceDeLEntropieDesProj}
$\\$If $T:(X,A)\to (Y,B)$ is measurable and if ${D_{\phi}}(P,Q)<\infty,$ then ${D_{\phi}}(P,Q)\geq {D_{\phi}}(PT^{-1},QT^{-1})$ with equality being reached when $T$ is surjective for $(P,Q)$.
\end{PROPIE}
\section{ Study of the sample}\label{truncSample}
Let $X_1$, $X_2$,..,$X_m$ be a sequence of independent random vectors with same density  $f$.
Let $Y_1$, $Y_2$,..,$Y_m$ be a sequence of independent random vectors with same density  $g$.
Then, the kernel estimators $f_m$, $g_m$, $f_{a,m}$ and $g_{a,m}$ of $f$, $g$, $f_{a}$ and $g_{a}$, for all $a\in\R^d_*$, almost surely and uniformly converge since we assume that the bandwidth $h_m$ of these estimators meets the following conditions (see \citet{MR0345296}): 

$(\mathcal Hyp)$: $h_m\searrow_m0$, $mh_m\nearrow_m\infty$, $mh_m/L(h_m^{-1})\to_m\infty$ and $L(h_m^{-1})/LLm\to_m\infty$,\\ with $L(u)=ln(u\vee e)$.\\
Let us consider 

$B_1(n,a)=\frac{1}{n}\Sigma_{i=1}^n\varphi'\{\frac{f_{a,n}(\transp aY_i)}{g_{a,n}(\transp aY_i) }\frac{g_n(Y_i)}{f_n(Y_i)}\}\frac{f_{a,n}(\transp aY_i)}{g_{a,n}(\transp aY_i) }$ and
$B_2(n,a)=\frac{1}{n}\Sigma_{i=1}^n\varphi^*\{\varphi'\{\frac{f_{a,n}(\transp aX_i)}{g_{a,n}(\transp aX_i) }\frac{g_n(X_i)}{f_n(X_i}\}\}.$\\
Our goal is to estimate the minimum of ${D_{\phi}}(g\frac{f_a}{g_a},f)$.
To do this, it is necessary for us to truncate our samples:\\
Let us consider now a positive sequence $\theta_m$ such that $\theta_m\to 0,\ y_m/\theta_n^2\to 0,$ where $y_m$ is the almost sure convergence rate of the kernel density estimator - $y_m=O_\PP(m^{-\frac{2}{4+d}})$, see lemma \ref{KernRate} - $y^{(1)}_m/\theta_m^2\to 0,$ where $y^{(1)}_m$ is defined by $$|\varphi(\frac{g_m(x)}{f_m(x)}\frac{f_{b,m}(\transp bx)}{g_{b,m}(\transp bx)})-\varphi(\frac{g(x)}{f(x)}\frac{f_b(\transp bx)}{g_b(\transp bx)})|\leq y^{(1)}_m$$ for all $b$ in $\R^d_*$ and all $x$ in $\R^d$, and finally  $\frac{y^{(2)}_m}{\theta_m^2}\to 0,$ where $y^{(2)}_n$ is defined by $$|\varphi'(\frac{g_m(x)}{f_m(x)}\frac{f_{b,m}(\transp bx)}{g_{b,m}(\transp bx)})-\varphi'(\frac{g(x)}{f(x)}\frac{f_b(\transp bx)}{g_b(\transp bx)})|\leq y^{(2)}_m$$ for all $b$ in $\R^d_*$ and all $x$ in $\R^d$.\\
We will generate $f_m$, $g_m$ and $g_{b,m}$ from the starting sample and we will select the $X_i$ and $Y_i$ vectors such that $f_m(X_i)\geq \theta_m$ and $g_{b,m}(\transp bY_i)\geq \theta_m$, for all $i$ and for all $b\in \R^d_*$. \\
The vectors meeting these conditions will be called $X_1,X_2,...,X_n$ and $Y_1,Y_2,...,Y_n$.\\
Consequently, the next proposition provides us with the condition required for us to derive our estimates:
\begin{proposition}\label{QuotientDonneLoi}
Using the notations introduced in \citet{MR2054155} and in section \ref{HypoF},
it holds 

$\lim_{n\to\infty}\sup_{a\in\R^d_*}|(B_1(n,a)-B_2(n,a))-{D_{\phi}}(g\frac{f_a}{g_a},f)|=0.$
\end{proposition}
\begin{remarque}\label{Scott}
With the Kullback-Leibler divergence, we can  take for $\theta_m$ the expression $m^{-\nu}$, with $0<\nu<\frac{1}{4+d}$.
\end{remarque}
\section{Hypotheses' discussion}\label{DiscussHyp}
\subsection{ Discussion of $(H1)$.}
\noindent Let us work with the Kullback-Leibler divergence and with $g$ and $a_1$. \\
For all $b\in\R^d_*$, we have
$\int\varphi^*(\varphi'(\frac{g(x)f_b(\transp bx)}{f(x)g_b(\transp bx)})) f(x)dx=\int (\frac{g(x)f_b(\transp bx)}{f(x)g_b(\transp bx)}-1)f(x)dx=0,$
since, for any $b$ in $\R^d_*$, the function $x\mapsto g(x)\frac{f_b(\transp bx)}{g_b(\transp bx)}$ is a density.
The complement of $\Theta^{D_{\phi}}$ in $\R^d_*$ is $\emptyset$ and then the supremum looked for in $\overline{\R}$ is $-\infty$.
We can therefore conclude.
It is interesting to note that we obtain the same verification with $f$, $g^{(k-1)}$ and $a_k$.
\subsection{ Discussion of $(H3)$.}
\noindent{\it This hypothesis consists in the following assumptions:\\
$\bullet$ We work with the Kullback-Leibler divergence, (0)\\
$\bullet$ We have $f(./\transp {a_1}x)=g(./\transp {a_1}x)$, i.e. $K(g\frac{f_1}{g_1},f)=0$ - we could also derive the same proof with $f$, $g^{(k-1)}$ and $a_k$ - (1)}\\
\noindent{\it Preliminary $(A)$:
Shows that
$A=\{(c,x)\in\R^d_*\backslash \{a_1\}\times R^d;\ \frac{f_{a_1}(\transp {a_1}x)}{g_{a_1}(\transp {a_1}x)}>\frac{f_{c}(\transp cx)}{g_{c}(\transp cx)},\ g(x)\frac{f_{c}(\transp cx)}{g_{c}(\transp cx)}> f(x)\}=\emptyset$
through a reductio ad absurdum, i.e. if we assume $A\not=\emptyset$.}\\
Thus, our hypothesis enables us to derive

$f(x)=f(./\transp {a_1}x)f_{a_1}(\transp {a_1}x)=g(./\transp {a_1}x)f_{a_1}(\transp {a_1}x)> g(./\transp {c}x)f_{c}(\transp {c}x)> f$\\
since $\frac{f_{a_1}(\transp {a_1}x)}{g_{a_1}(\transp {a_1}x)}\geq\frac{f_{c}(\transp cx)}{g_{c}(\transp cx)}$ implies $g(./\transp {a_1}x)f_{a_1}(\transp {a_1}x)=g(x)\frac{f_{a_1}(\transp {a_1}x)}{g_{a_1}(\transp {a_1}x)}\geq g(x)\frac{f_{c}(\transp cx)}{g_{c}(\transp cx)}=g(./\transp {c}x)f_{c}(\transp {c}x)$, i.e. $f>f$. We can therefore conclude.
$\\${\it Preliminary $(B)$:
Shows that $B=\{(c,x)\in\R^d_*\backslash \{a_1\}\times R^d;\ \frac{f_{a_1}(\transp {a_1}x)}{g_{a_1}(\transp {a_1}x)}<\frac{f_{c}(\transp cx)}{g_{c}(\transp cx)},\ g(x)\frac{f_{c}(\transp cx)}{g_{c}(\transp cx)}< f(x)\}=\emptyset$
through a reductio ad absurdum, i.e. if we assume $B\not=\emptyset$.}\\
Thus, our hypothesis enables us to derive

$f(x)=f(./\transp {a_1}x)f_{a_1}(\transp {a_1}x)=g(./\transp {a_1}x)f_{a_1}(\transp {a_1}x)< g(./\transp {c}x)f_{c}(\transp {c}x)< f$\\
We can therefore conclude as above.
$\\$Let us now verify  $(H3)$:\\
We have $P M(c,a_1)- P M(c,a)=\int ln(\frac{g(x)f_{c}(\transp cx)}{g_{c}(\transp cx)f(x)})\{\frac{f_{a_1}(\transp {a_1}x)}{g_{a_1}(\transp {a_1}x)}-\frac{f_{c}(\transp cx)}{g_{c}(\transp cx)}\}g(x)dx.$
Moreover, the logarithm $ln$ is negative on 
$\{x\in\R^d_*;\ \frac{g(x)f_{c}(\transp cx)}{g_{c}(\transp cx)f(x)}<1\}$ and is positive on $\{x\in\R^d_*;\ \frac{g(x)f_{c}(\transp cx)}{g_{c}(\transp cx)f(x)}\geq1\}$.\\
Thus, the preliminary studies $(A)$ and $(B)$ show that $ln(\frac{g(x)f_{c}(\transp cx)}{g_{c}(\transp cx)f(x)})$ and $\{\frac{f_{a_1}(\transp {a_1}x)}{g_{a_1}(\transp {a_1}x)}-\frac{f_{c}(\transp cx)}{g_{c}(\transp cx)}\}$  always present a negative product. We can therefore conclude, since $(c,a)\mapsto P M(c,a_1)- P M(c,a)$ is not null for all $c$ and for all $a$ - with $a\not=a_1$.

\section{On Huber's algorithms}\label{PP-old}
In the present appendix, let us now first present the projection pursuit methodologies introduced by \citet{MR790553}. Secondly, we will show that our method encompasses Hubers'.

\subsection{ Huber's analytic approach}
Let $f$ be a density on $\R^d$. We define an instrumental density $g$ with the same mean and variance as  $f$.
Huber's methodology requires us to start with performing the $K(f,g)=0$ test - with $K$ being the Kullback-Leibler divergence. Should this test turn out to be positive, then $f=g$ and the algorithm stops. If the test were not to be verified, the first step of Huber's algorithm would amount to defining a vector $a_1$ and a density $f^{(1)}$  by 
\begin{equation}\label{DefSequMethodH1}
a_1\ =\ arg\inf_{a\in\R^d_*}\ K(f\frac{g_a}{f_a},g)\text{ and }f^{(1)}=f\frac{g_{a_1}}{f_{a_1}}
\end{equation}
where $\R^d_*$ is the set of non null vectors of $\R^d$ and $f_a$ (resp. $g_a$) stands for the density of $\transp aX$ (resp. $\transp aY$) when $f$ (resp. $g$) is the density of $X$ (resp. $Y$). More exactly, this results from the maximisation of $a\mapsto K(f_a,g_a)$ since $K(f,g)=K(f_a,g_a)+K(f\frac{g_a}{f_a},g)$ and it is assumed that $K(f,g)$ is finite.
In a second step, Huber replaces  $f$ with $f^{(1)}$ and goes through the first step again.\\
By iterating this process, Huber thus obtains a sequence $(a_1,a_2,...)$ of vectors of $\R^d_*$ and a sequence of densities $f^{(i)}$.
\begin{remarque}
This algorithm stops when the Kullback-Leibler divergence equals zero or when it reaches the $d^{th}$ iteration. We then obtain an approximation of $f$ from $g$ :\\
When there exists an integer $j$ such that $K(f^{(j)},g)=0$ with $j\leq d$, he obtains $f^{(j)}=g$, i.e. $f=g\Pi_{i=1}^j\frac{f^{(i-1)}_{a_i}} {g_{a_i}}$ since by induction $f^{(j)}=f\Pi_{i=1}^j\frac{g_{a_i}}{f^{(i-1)}_{a_i}}$. Similarly, when, for all $j$, Huber gets $K(f^{(j)},g)>0$ with $j\leq d$, he assumes $g=f^{(d)}$ in order to derive $f=g\Pi_{i=1}^d\frac{f^{(i-1)}_{a_i}} {g_{a_i}}$. \\
Finally, he obtains  $K(f^{(0)},g)\geq K(f^{(1)},g)\geq.....\geq 0$ with $f^{(0)}=f$.
\end{remarque}
\subsection{ Huber's synthetic approach}
Keeping the notations of the above section, we start with performing the $K(f,g)=0$ test; should this test turn out to be positive, then $f=g$ and the algorithm stops, otherwise, the first step of his algorithm would consist in defining a vector $a_1$ and a density $g^{(1)}$  by 
\begin{equation}\label{DefSequMethodH2}
a_1\ =\ arg\inf_{a\in\R^d_*}\ K(f,g\frac{f_a}{g_a})\text{ and }g^{(1)}=g\frac{f_{a_1}}{g_{a_1}}
\end{equation}
More exactly, this optimisation results from the maximisation of $a\mapsto K(f_a,g_a)$ since $K(f,g)=K(f_a,g_a)+K(f,g\frac{f_a}{g_a})$ and it is assumed that $K(f,g)$ is finite.
In a second step, Huber replaces  $g$ with $g^{(1)}$ and goes through the first step again.
By iterating this process, Huber thus obtains a sequence $(a_1,a_2,...)$ of vectors of $\R^d_*$ and a sequence of densities $g^{(i)}$.
\begin{remarque}
First, in a similar manner to the analytic approach, this methodology enables us to approximate $f$ from $g$:\\
To obtain an approximation of $f$, Huber either stops his algorithm when the Kullback-Leibler divergence equals zero, i.e. $K(f,g^{(j)})=0$ implies  $g^{(j)}=f$ with $j\leq d$, or when his  algorithm reaches the $d^{th}$ iteration, i.e. he approximates $f$ with $g^{(d)}$.\\
Second, he gets $K(f,g^{(0)})\geq K(f,g^{(1)})\geq.....\geq 0$ with $g^{(0)}=g$.
\end{remarque}
\subsection{ The first co-vector of $f$ simultaneously optimizes four problems}\label{a1solve4}
Let us first study Huber's analytic approach.\\
Let $\cR'$ be the class of all positive functions $r$ defined on $\R$ and such that $f(x)r^{-1}(\transp ax)$ is a density on $\R^d$  for all $a$ belonging to $\R^d_*$. The following proposition shows that there exists a vector $a$ such that  $\frac{f_a}{g_a}$ minimizes $K(fr^{-1},g)$ in $r$:
\begin{proposition}[Analytic Approach]\label{lemmeHuber0prop}
There exists a vector $a$ belonging to $\R^d_*$ such that 

$arg\min_{r\in\cR'}K(fr^{-1},g)=\frac{f_a}{g_a},\text{ and } r(\transp ax)=\frac{f_a(\transp ax)}{g_a(\transp ax)}$
and $  K(f,g)=  K(f_a,g_a)+  K(f\frac{g_a}{f_a},g).$
\end{proposition}
\noindent Let us also study Huber's synthetic approach:\\
Let $\cR$ be the class of all positive functions $r$ defined on $\R$ and such that $g(x)r(\transp ax)$ is a density on $\R^d$ for all $a$ belonging to $\R^d_*$. The following proposition shows that there exists a vector $a$ such that  $\frac{f_a}{g_a}$ minimizes $K(gr,f)$ in $r$:
\begin{proposition}[Synthetic Approach]\label{lemmeHuberprop}
There exists a vector $a$ belonging to $\R^d_*$ such that

$arg\min_{r\in\cR}K(f,gr)=\frac{f_a}{g_a},\text{ and } r(\transp ax)=\frac{f_a(\transp ax)}{g_a(\transp ax)}$
and $  K(f,g)=  K(f_a,g_a)+  K(f,g\frac{f_a}{g_a}).$
\end{proposition}
To recapitulate, the choice of $r=\frac{f_a}{g_a}$ enables us to simultaneously solve the following three optimisation problems, for $a\in\R^d_*$:

$\text{First, find $a$ such that }a\ =\ {\text{\itshape arginf}}_{a\in\R^d_*}\ K(f\frac{g_a}{f_a},g)$ 
- analytic approach -

$\text{Second, find $a$ such that }a\ =\ {\text{\itshape arginf}}_{a\in\R^d_*}\ K(f,g\frac{f_a}{g_a})$ - synthetic approach -

$\text{Third, find $a$ such that }a\ =\ {\text{\itshape arginf}}_{a\in\R^d_*}\ D_\phi(g\frac{f_a}{g_a},f)$ - our method.\\
We can therefore state that the methodology we introduced in the present article encompasses Hubers'.
\section{Proofs}
\noindent{\bf Proof of propositions \ref{lemmeHuber0prop} and \ref{lemmeHuberprop}.}
Let us first study proposition \ref{lemmeHuberprop}.\\
Without loss of generality, we will prove this proposition with $x_1$ in lieu of $\transp aX$.\\
Let us define $g^*=gr$. We remark that $g$ and $g^*$ present the same density conditionally to $x_1$.
Indeed,\\ $g^*_1(x_1)= \int g^*(x)dx_2...dx_d= \int r(x_1)g(x)dx_2...dx_d= r(x_1)\int g(x)dx_2...dx_d= r(x_1)g_1(x_1)$.\\
Thus, we can demonstrate this proposition.\\
We have $g(.|x_1)=\frac{g(x_1,..., x_n)}{g_1(x_1)}$ and $g_1(x_1)r(x_1)$ is the marginal density of $g^*.$
Hence,\\ $\int g^*dx=\int g_1(x_1)r(x_1)g(.|x_1)dx=\int g_1(x_1)\frac{f_1(x_1)}{g_1(x_1)}(\int g(.|x_1)dx_2..dx_d)dx_1=\int f_1(x_1)dx_1=1$ and since $g^*$ is positive, then $g^*$ is a density.
Moreover,
\begin{eqnarray}
  K(f,g^*)&=& \int f\{ln(f)-ln(g^*)\}dx,\label{dernière-1}\\
        &=& \int f\{ln(f(.|x_1))-ln(g^*(.|x_1))+ln(f_1(x_1))-ln(g_1(x_1)r(x_1))\}dx,\nonumber\\
        &=& \int f\{ln(f(.|x_1))-ln(g(.|x_1)) +ln(f_1(x_1))-ln(g_1(x_1)r(x_1))\}dx,\label{dernière}
\end{eqnarray}
as $g^*(.|x_1)=g(.|x_1)$. Since the minimum of this last equation (\ref{dernière}) is reached through the minimization of $\int f\{ln(f_1(x_1))-ln(g_1(x_1)r(x_1))\}dx=  K(f_1,g_1r)$,
then property  \ref{Phimini} necessarily implies that $f_1=g_1r$, hence $r=f_1/g_1$.\\
Finally, we have $K(f,g)-  K(f,g^*)= \int f\{ln(f_1(x_1))-ln(g_1(x_1))\}dx=K(f_1,g_1),$
which completes the demonstration of proposition \ref{lemmeHuberprop}.\\
Similarly, if we replace $f^*=fr^{-1}$ with $f$ and $g$ with $g^*$, we obtain the proof of proposition \ref{lemmeHuber0prop}.\hfill$\Box$\\
\noindent{\bf Proof of proposition \ref{lemmeHuberModifprop}.}
The demonstration is also very similar to the one for proposition \ref{lemmeHuberprop}, save for the fact we now base our reasoning at row \ref{dernière-1} on 
$K(g^*,f)= \int g^*\{ln(f)-ln(g^*)\}dx$ instead of $K(f,g^*)= \int f\{ln(f)-ln(g^*)\}dx$.\hfill$\Box$\\
\noindent{\bf Proof of lemma \ref{H3}.}
\begin{lemme}\label{H3}
We have $\Theta=\{b\in\Theta\ |\ \ \int(\frac{g(x)}{f(x)}\frac{f_b(\transp bx)}{g_b(\transp bx)}-1)f(x)dx<\infty\}$.
\end{lemme}
We get the result since 
$\int\ (\frac{g(x)f_b(\transp bx)}{f(x)g_b(\transp bx)}-1)f(x)\ dx=\int\ (\frac{g(x)f_b(\transp bx)}{g_b(\transp bx)}-f(x))\ dx=0$.\hfill$\Box$\\
\noindent{\bf Proof of proposition \ref{pConv2}.}
\begin{proposition}\label{pConv2}
In the case where $f$ is known and keeping the notations introduced in section \ref{HypoF}, as well as
assuming $(H1)$ to $(H3)$ hold, then both $\sup_{a\in\Theta}\|\check c_n(a)-a_k\|$ and $\check \gamma_n $ tend to $a_k$ a.s.
\end{proposition}
In the same manner as in Proposition 3.4 of \citet{MR2054155}, we prove this proposition through lemma \ref{H3}.\hfill$\Box$\\

\noindent{\bf Proof of proposition \ref{KernelpConv2}.}
Proposition \ref{KernelpConv2} comes immediately from proposition \ref{QuotientDonneLoi} page \pageref{QuotientDonneLoi} and lemma \ref{pConv2} page \pageref{pConv2}.\hfill$\Box$\\
\noindent{\bf Proof of theorem \ref{KernelKRessultatPricipal}.}
We prove this theorem by induction.
First, by the very definition of the kernel estimator $\check g^{(0)}_n=g_n$ converges towards $g$. Moreover, the continuity of $a\mapsto f_{a,n}$ and $a\mapsto g_{a,n}$ and proposition \ref{KernelpConv2} imply that $\check g^{(1)}_n=\check g^{(0)}_n\frac{f_{a,n}}{\check g^{(0)}_{a,n}}$ converges towards $g^{(1)}$. 
Finally, since, for any $k$, $\check g^{(k)}_n=\check g^{(k-1)}_n\frac{f_{\check a_k,n}}{\check g^{(k-1)}_{\check a_k,n}}$, we conclude similarly as for $\check g^{(1)}_n$.\hfill$\Box$\\
\noindent{\bf Proof of lemma \ref{ChangBasis}.}
\begin{lemme} \label{ChangBasis}
We have $g(./\transp{a_{1}}x,...,\transp{a_{j}}x)=n(\transp{a_{j+1}}x,...,\transp{a_{d}}x)=f(./\transp{a_{1}}x,...,\transp{a_{j}}x)$.
\end{lemme}
\noindent Putting $A=(a_1,..,a_d)$, let us determine $f$ in basis $A$.
Let us first study the function defined by $\psi:\R^d\to\R^d$, $x\mapsto(\transp{a_1}x,..,\transp{a_d}x).$
We can immediately say that $\psi$ is continuous and since $A$ is a basis, its bijectivity is obvious.
Moreover, let us study its Jacobian.\\
By definition, it is $J_\psi(x_1,\dotsc,x_d)=
\begin{vmatrix}
\displaystyle\frac{\partial\psi_1}{\partial x_1} & \dotsb & \displaystyle\frac{\partial\psi_1}{\partial x_d}\\
\dotsb & \dotsb & \dotsb\\
\displaystyle\frac{\partial\psi_d}{\partial x_1} & \dotsb & 
\displaystyle\frac{\partial\psi_d}{\partial x_d}
\end{vmatrix}=
\begin{vmatrix}
\displaystyle a_{1,1} & \dotsb & \displaystyle a_{1,d}\\
\dotsb & \dotsb & \dotsb\\
\displaystyle a_{d,1} & \dotsb & 
\displaystyle a_{d,d}
\end{vmatrix}=|A|\not=0$ since $A$ is a basis. We can therefore infer : 
$\forall x \in\R^d,\ \exists! y\in\R^d\text{ such that }f(x)=|A|^{-1}\Psi(y),$
i.e. $\Psi$ (resp. $y$) is the expression of $f$ (resp of $x$) in  basis $A$, namely
$\Psi(y)=\tilde n(y_{j+1},...,y_{d})\tilde h(y_{1},...,y_{j})$, with $\tilde n$ and $\tilde h$ being the expressions of $n$ and $h$ in basis $A$. 
Consequently, our results in the case where the family $\{a_j\}_{1\leq j\leq d}$ is the canonical basis of $\R^d$, still hold for $\Psi$ in  basis $A$ - see section \ref{modelSection}. And then, if $\tilde g$ is the expression of $g$ in basis $A$, we have
$\tilde g(./y_1,...,y_{j})= \tilde n(y_{j+1},...,y_{d})=\Psi(./y_{1},...,y_{j})$, i.e. $g(./\transp{a_{1}}x,...,\transp{a_{j}}x)=n(\transp{a_{j+1}}x,...,\transp{a_{d}}x)=f(./\transp{a_{1}}x,...,\transp{a_{j}}x)$.\hfill$\Box$\\

\noindent{\bf Proof of lemma \ref{KernRate}.}
\begin{lemme}\label{KernRate}
For any continuous density $f$, we have
$y_m=|f_m(x)-f(x)|=O_\PP(m^{-\frac{2}{4+d}})$.
\end{lemme}
Defining $b_m(x)$ as  $b_m(x)=|E(f_m(x))-f(x)|$, we have $y_m\leq |f_m(x)-E(f_m(x))|+b_m(x)$. Moreover, 
from page 150 of \citet{MR1191168}, we derive that $b_m(x)=O_\PP(\Sigma_{j=1}^dh_j^2)$ where $h_j=O_\PP(m^{-\frac{1}{4+d}})$. Then, we obtain $b_m(x)=O_\PP(m^{-\frac{2}{4+d}})$. Finally, since the central limit theorem rate  is $O_\PP(m^{-\frac{1}{2}})$, we infer that $y_m\leq O_\PP(m^{-\frac{1}{2}})+O_\PP(m^{-\frac{2}{4+d}})=O_\PP(m^{-\frac{2}{4+d}})$.\hfill$\Box$\\
\noindent{\bf Proof of lemma \ref{aleph}.} 
\begin{lemme}\label{aleph}
Let $f$ be an absolutely continuous density, then, for all sequences $(a_n)$ tending  to $a$ in $\R^d_*$, sequence $f_{a_n}$ uniformly converges towards $f_a$.
\end{lemme}
\begin{proof}
For all $a$ in $\R^d_*$, let $F_a$ be the cumulative distribution function of $\transp aX$ and $\psi_a$ be a complex function defined by $\psi_a(u,v)=F_a(\cR e(u+iv))+iF_a(\cR e(v+iu))$, for all $u$ and $v$ in $\R$.\\First, the function $\psi_a(u,v)$ is an analytic function, because $x\mapsto f_a(\transp ax)$ is continuous and as a result of the corollary of Dini's second theorem - according to which 
{\it "A sequence of cumulative distribution functions, which pointwise converges on $\R$ towards a continuous cumulative distribution function $F$ on $\R$, uniformly converges towards $F$ on $\R$"}-
we deduct that, for all sequences $(a_n)$ converging towards $a$, $\psi_{a_n}$ uniformly converges towards $\psi_a$.
Finally, the Weierstrass theorem, (see proposal $(10.1)$ page 220 of the "Calcul infinit\'esimal" book of Jean Dieudonn\'e), implies that all sequences $\psi_{a,n}'$ uniformly converge towards $\psi_a'$, for all $a_n$ tending to $a$. We can therefore conclude.
\end{proof}

\noindent{\bf Proof of lemma \ref{compacité-1}.} 
By definition of the closure of a set, we have
\begin{lemme} \label{compacité-1}
The set $\Gamma_c$ is closed in $L^1$ for the topology of the uniform convergence.
\end{lemme}
\noindent{\bf Proof of lemma \ref{compacité}.} 
Since $K$ is greater than the $L^1$ distance, we have
\begin{lemme} \label{compacité}
For all $c>0$, we have $\Gamma_c\subset \overline B_{L^1}(f,c),$ where $B_{L^1}(f,c)=\{p\in L^1;\|f-p\|_1\leq c\}$.
\end{lemme}
\noindent{\bf Proof of lemma \ref{compacité+1}.} 
The definition of the closure of a set and lemma \ref{aleph} (see page \pageref{aleph}) imply :
\begin{lemme} \label{compacité+1}
 $G$ is closed in $L^1$ for the topology of the uniform convergence.
\end{lemme}
\noindent{\bf Proof of lemma \ref{toattain}.}
\begin{lemme}\label{toattain}
$\inf_{a\in\R^d_*}  {D_{\phi}}(g^*,f)$ is reached when the $\phi$-divergence is greater than the $L^1$ distance as well as the $L^2$ distance.
\end{lemme}
\begin{proof}
Indeed, let $G$ be $\{g\frac{f_a}{g_a};\ a\in\R^d_*\}$ and $\Gamma_c$ be $\Gamma_c=\{p;\   K(p,f)\leq c\}$ for all $c$>0. From lemmas \ref{compacité-1}, \ref{compacité} and \ref{compacité+1} (see page \pageref{compacité}), we get $\Gamma_c\cap G$ is a compact for the topology of the uniform convergence, if $\Gamma_c\cap G$ is not empty.
Hence, and since  property \ref{K-SCI} (see page \pageref{K-SCI}) implies that  $Q\mapsto   {D_{\phi}}(Q,P)$ is lower semi-continuous in $L^1$ for the topology of the uniform convergence, then the infimum is reached in $L^1$.
(Taking for example  $c={D_{\phi}}(g,f),$ $\Omega$  is necessarily not empty because we always have ${D_{\phi}}(g\frac{f_a}{g_a},f)\leq {D_{\phi}}(g,f)$).
Moreover, when the $\phi-$divergence is greater than the $L^2$ distance, the very definition of the $L^2$ space enables us to provide the same proof as for the $L^1$ distance.
\end{proof}

\noindent{\bf Proof of lemma \ref{TrucBidule}.}
\begin {lemme}\label{TrucBidule}
For any $p\leq d$, we have $g^{(p-1)}_{a_p}=g_{a_p}$.
\end{lemme}
Assuming, without any loss of generality, that the $a_i$, $i=1,..,p$, are the vectors of the canonical basis, since 
$g^{(p-1)}(x)=g(x)\frac{f_1(x_1)}{g_1(x_1)}\frac{f_2(x_2)}{g_2(x_2)}...\frac{f_{p-1}(x_{p-1})}{g_{p-1}(x_{p-1})}$ we derive immediately that $g^{(p-1)}_{p}=g_{p}$. We note that it is sufficient to operate a change in basis on the $a_i$ to obtain the general case since $A=(a_i)$ is a basis - see lemma \ref{imFree}.\hfill$\Box$\\
\noindent{\bf Proof of lemma \ref{imFree}.}
\begin {lemme}\label{imFree}
If there exists $p$, $p\leq d$, such that ${D_{\phi}}(g^{(p)},f)=0$, then the family  of $(a_i)_{i=1,..,p}$ is free and is orthogonal.
\end{lemme}
Without any loss of generality, let us assume that $p=2$ and that the $a_i$ are the vectors of the canonical basis. Using a reductio ad absurdum based on the hypotheses  $a_1=(1,0,...,0)$ and $a_2=(\alpha,0,...,0)$, where  $\alpha\in\R$, we get $g^{(1)}(x)=g(x_2,..,x_d/x_1)f_1(x_1)$ and $f=g^{(2)}(x)=g(x_2,..,x_d/x_1)f_1(x_1)\frac{f_{\alpha a_1}(\alpha x_1)}{[g^{(1)}]_{\alpha a_1}(\alpha x_1)}$. Hence $f(x_2,..,x_d/x_1)=g(x_2,..,x_d/x_1)\frac
{f_{\alpha a_1}(\alpha x_1)}
{[g^{(1)}]_{\alpha a_1}(\alpha x_1)}.$
It consequently implies that $f_{\alpha a_1}(\alpha x_1)=[g^{(1)}]_{\alpha a_1}(\alpha x_1)$ since
 $1=\int f(x_2,..,x_d/x_1)dx_2...dx_d=\int g(x_2,..,x_d/x_1)dx_2...dx_d\frac
{f_{\alpha a_1}(\alpha x_1)}
{[g^{(1)}]_{\alpha a_1}(\alpha x_1)}=\frac
{f_{\alpha a_1}(\alpha x_1)}
{[g^{(1)}]_{\alpha a_1}(\alpha x_1)}$.\\
Therefore, $g^{(2)}=g^{(1)}$, i.e. $p=1$ which leads to a contradiction. Hence, the family is free.\\
Moreover, using a reductio ad absurdum, we get the orthogonality. Indeed, we have \\
$\int f(x)dx=1\not=+\infty=\int n(\transp{a_{j+1}}x,...,\transp{a_{d}}x)h(\transp{a_{1}}x,...,\transp{a_{j}}x)dx$. The use of the same argument as in the proof of lemma \ref{Base}, enables us to infer the orthogonality of $(a_i)_{i=1,..,p}$.\hfill$\Box$\\
\noindent{\bf Proof of lemma \ref{Base}.}
\begin{lemme} \label{Base}
Should there exist a family $(a_i)_{i=1...d}$ such that 
$f(x)=n(\transp{a_{j+1}}x,...,\transp{a_{d}}x)h(\transp{a_{1}}x,...,\transp{a_{j}}x),$
with $j<d$, with $f$, $n$ and $h$ being densities, then this family is an orthogonal basis of $\R^d$.
\end{lemme}
\noindent Using a reductio ad absurdum, we have $\int f(x)dx=1\not=+\infty=\int n(\transp{a_{j+1}}x,...,\transp{a_{d}}x)h(\transp{a_{1}}x,...,\transp{a_{j}}x)dx$. We can therefore conclude.\hfill$\Box$\\
\noindent{\bf Proof of proposition \ref{PropCop}.}\\
Through lemma \ref{imFree}, we can consequently infer  that $(a_1,...,a_d)$ is a basis of $\R^d.$
Let us now write $f$ in the $A$ system.
Let us first study the function defined by $\psi\ :\ \R^d\to\R^d,\ x\mapsto(\transp{a_1}x,..,\transp{a_d}x).$
We can say $\psi$ is continuous and since $A$ is a basis, its bijectivity is obvious.
Let us also study its Jacobian. By definition, it is
$
J_\psi(x_1,\dotsc,x_k)=
\begin{vmatrix}
\displaystyle\frac{\partial\psi_1}{\partial x_1} & \dotsb & \displaystyle\frac{\partial\psi_1}{\partial x_k}\\
\dotsb & \dotsb & \dotsb\\
\displaystyle\frac{\partial\psi_k}{\partial x_1} & \dotsb & 
\displaystyle\frac{\partial\psi_k}{\partial x_k}
\end{vmatrix}=
\begin{vmatrix}
\displaystyle a_{1,1} & \dotsb & \displaystyle a_{1,d}\\
\dotsb & \dotsb & \dotsb\\
\displaystyle a_{d,1} & \dotsb & 
\displaystyle a_{d,d}
\end{vmatrix}=|A|\not=0
$
since $A$ is a basis. Thus, we can infer that, in basis $A$, the writing of $f$ (resp. $x$) exists and is unique. Defining $\Psi$ (resp. $y$) as this new form of $f$ (resp. $x$), we have $f(x)=|A|^{-1}\Psi(y)$ (resp. $(\transp{a_1}x,..,\transp{a_d}x)=(y_1,...,y_d)$). Similarly, let us define $\tilde \Psi$ (resp. $\tilde \Psi^{(i)}$) as being the form of $g$ (resp. $g^{(i)}$) in basis $A$, we also have $\tilde \Psi(y)=|A|g(x)$ \\(resp.  $\tilde \Psi^{(i)}(y)=|A|g^{(i)}(x)$).\\
Now, through a finite induction in $i$, $1\leq i\leq d$, let us demonstrate the following property

$\PP(i)="\tilde \Psi^{(i)}(y)=\tilde \Psi(y_{i+1},...,y_d/y_1,...,y_i)\Psi_1(y_1)\Psi_2(y_2)...\Psi_i(y_i)"$ \\
Initialisation :

For $i=1$. The above notations lead us to
$\tilde \Psi^{(1)}(y)=\tilde \Psi(y)\frac{\Psi_1(y_1)}{\tilde\Psi_1(y_1)}$, 
since $y_1=\transp{a_1}x$ through the change in variables, i.e.
$\tilde \Psi^{(1)}(y)=\tilde \Psi(y_2,...,y_d/y_1)\Psi_1(y_1),$ by the very definition of conditional density.
Hence, $\PP(1)$ holds true.

For $i=2$. Since $\PP(1)$ is true, we can write\\
$\tilde \Psi^{(2)}(y)=\tilde \Psi^{(1)}(y)\frac{\Psi_2(y_2)}{\tilde\Psi^{(1)}_2(y_2)}\text{, by  definition of $\tilde \Psi^{(2)}(y)$,}$

$=\tilde \Psi(y_2,...,y_d/y_1)\Psi_1(y_1)\frac{\Psi_2(y_2)}{\tilde\Psi^{(1)}_2(y_2)}\text{, since $\PP(1)$ is true,}$

$=\tilde \Psi(y_2,...,y_d/y_1)\Psi_1(y_1)\frac{\Psi_2(y_2)}{\tilde\Psi_2(y_2)}\text{, since $\tilde\Psi^{(1)}_2(y_2)=\tilde\Psi_2(y_2)$ ,}$

$=\tilde \Psi(y_3,...,y_d/y_1,y_2)\Psi_1(y_1)\Psi_2(y_2),$\\
by the very definition of conditional density. Thus, $\PP(2)$ holds true.\\
Going from $i-1$ to $i$ $(i\leq p)$:\\
Let us assume $\PP(i-1)$ is true, we can then show that $\PP(i)$.\\
$\tilde \Psi^{(i)}(y)=\tilde \Psi^{(i-1)}(y)\frac{\Psi_i(y_i)}{\tilde\Psi^{(i-1)}_i(y_i)}\text{, by definition of $\tilde \Psi^{(i)}(y)$,}$

$=\tilde \Psi(y_i,...,y_d/y_1,..,y_{i-1})\Psi_1(y_1)...\Psi_{i-1}(y_{i-1})\frac{\Psi_i(y_i)}{\tilde\Psi^{(i-1)}_i(y_i)}$,$\text{ since $\PP(i-1)$ is true,}$

$=\tilde \Psi(y_i,...,y_d/y_1,...,y_{i-1})
\Psi_1(y_1)...\Psi_{i-1}(y_{i-1})
\frac{\Psi_i(y_i)}{\tilde\Psi_i(y_i)},\text{ since $\tilde\Psi^{(i-1)}_i(y_i)=\tilde\Psi_i(y_i)$ ,}$

$=\tilde \Psi(y_{i+1},...,y_d/y_1,...,y_i)\Psi_1(y_1)\Psi_2(y_2)...\Psi_i(y_i),\text{ by the very definition of conditional density.}$\\
Thus, $\PP(i)$ is true.\\
Conclusion :\\
The induction principle enables us to infer that $\PP(i)$ holds true for $1\leq i\leq d$.\\

\noindent At present, since ${D_{\phi}}(\tilde \Psi^{(d)},\Psi)=0$, the above entails that $\tilde \Psi^{(d)}(y)=\Psi(y)$,
i.e. 

$\frac{\tilde \Psi(y)} {\tilde \Psi_1(y_1)\tilde \Psi_2(y_2)...\tilde \Psi_d(y_d)}=\frac{\Psi(y)}{\Psi_1(y_1) \Psi_2(y_2)...\Psi_d(y_d)}.$\\
We finally obtain $\frac{\dr^d}{\dr y_1...\dr y_d}\tilde C=\frac{\dr^d}{\dr y_1...\dr y_d}C$, 
where $C$ and $\tilde C$ are the respective copulas of $\Psi$ and $\tilde \Psi.$\\
Let us remark that, if the $(a_i)$ are the canonical basis of $\R^d$, we have

$\frac{\dr^d}{\dr x_1...\dr x_d}C_f=\frac{\dr^d}{\dr x_1...\dr x_d}C_{g}$,\\
where $C_f$ and $C_{g}$ are the respective copulas of $f$ and $g.$\hfill$\Box$
$\\$\noindent{\bf Proof of proposition  \ref{QuotientDonneLoi}.}\\
Let us first note that we will prove this proposition for $k\geq2$, i.e. in the case where $g^{(k-1)}$ is not known. The initial case using the known density $g^{(0)}=g$, will be an immediate consequence of the above. \\
Moreover, going forward, to be more legible, we will use $g$ (resp. $g_n$) in lieu of $g^{(k-1)}$ (resp. $g^{(k-1)}_n$).\\
We can therefore remark that we have $f(X_i)\geq \theta_n-y_n$, $g(Y_i)\geq \theta_n-y_n$ and $g_b(\transp bY_i)\geq \theta_n-y_n$, for all $i$ and for all $b\in \R^d_*$, thanks to the uniform convergence of the kernel estimators.
Indeed, we have $f(X_i)=f(X_i)-f_n(X_i)+f_n(X_i)\geq-y_n+f_n(X_i)$, by definition of $y_n$, and then  $f(X_i)\geq-y_n+\theta_n$, by hypothesis on $f_n(X_i)$. This is also true for $g_n$ and $g_{b,n}$.\\
This entails
$\sup_{b\in\R^d_*}|\frac{1}{n}\Sigma_{i=1}^{n}\varphi'(\frac{f_{a,n}(\transp aY_i)}{g_{a,n}(\transp aY_i) }\frac{g_n(Y_i)}{f_n(Y_i)}).\frac{f_{a,n}(\transp aY_i)}{g_{a,n}(\transp aY_i) }-\int \varphi'(\frac{g(x)}{f(x)}\frac{f_b(\transp bx)}{g_b(\transp bx)}) g(x)\frac{f_a(\transp ax)}{g_a(\transp ax)} dx|\to0\ a.s.$\\
Indeed, let us remark that

$|\frac{1}{n}\Sigma_{i=1}^{n}\{\varphi'\{\frac{f_{a,n}(\transp aY_i)}{g_{a,n}(\transp aY_i) }\frac{g_n(Y_i)}{f_n(Y_i)}\}\frac{f_{a,n}(\transp aY_i)}{g_{a,n}(\transp aY_i) }\}-\int\ \varphi'(\frac{g(x)}{f(x)}\frac{f_b(\transp bx)}{g_b(\transp bx)})\ g(x)\frac{f_a(\transp ax)}{g_a(\transp ax)}\ dx|$

$=|\frac{1}{n}\Sigma_{i=1}^{n}\varphi'\{\frac{f_{a,n}(\transp aY_i)}{g_{a,n}(\transp aY_i) }\frac{g_n(Y_i)}{f_n(Y_i)}\}\frac{f_{a,n}(\transp aY_i)}{g_{a,n}(\transp aY_i) }-\frac{1}{n}\Sigma_{i=1}^{n}\varphi'\{\frac{f_{a}(\transp aY_i)}{g_{a}(\transp aY_i) }\frac{g(Y_i)}{f(Y_i)}\}\frac{f_{a}(\transp aY_i)}{g_{a}(\transp aY_i) }$

$\ \ \ \ \ $ $+\frac{1}{n}\Sigma_{i=1}^{n}\varphi'\{\frac{f_{a}(\transp aY_i)}{g_{a}(\transp aY_i) }\frac{g(Y_i)}{f(Y_i)}\}\frac{f_{a}(\transp aY_i)}{g_{a}(\transp aY_i) }-\int\ \varphi'(\frac{g(x)}{f(x)}\frac{f_b(\transp bx)}{g_b(\transp bx)})\ g(x)\frac{f_a(\transp ax)}{g_a(\transp ax)}\ dx|$

$\leq|\frac{1}{n}\Sigma_{i=1}^{n}\varphi'\{\frac{f_{a,n}(\transp aY_i)}{g_{a,n}(\transp aY_i) }\frac{g_n(Y_i)}{f_n(Y_i)}\}\frac{f_{a,n}(\transp aY_i)}{g_{a,n}(\transp aY_i) }-\frac{1}{n}\Sigma_{i=1}^{n}\varphi'\{\frac{f_{a}(\transp aY_i)}{g_{a}(\transp aY_i) }\frac{g(Y_i)}{f(Y_i)}\}\frac{f_{a}(\transp aY_i)}{g_{a}(\transp aY_i) }|$

$\ \ \ \ \ $ $+|\frac{1}{n}\Sigma_{i=1}^{n}\varphi'\{\frac{f_{a}(\transp aY_i)}{g_{a}(\transp aY_i) }\frac{g(Y_i)}{f(Y_i)}\}\frac{f_{a}(\transp aY_i)}{g_{a}(\transp aY_i) }-\int\ \varphi'(\frac{g(x)}{f(x)}\frac{f_b(\transp bx)}{g_b(\transp bx)})\ g(x)\frac{f_a(\transp ax)}{g_a(\transp ax)}\ dx|$\\
Moreover, since $\int|\varphi'(\frac{g(x)}{f(x)}\frac{f_b(\transp bx)}{g_b(\transp bx)})\ g(x)\frac{f_a(\transp ax)}{g_a(\transp ax)}|dx< \infty$, as implied by lemma \ref{ExitenceDeLEntropieDesProj}, and since we assumed $g$ such that ${D_{\phi}}(g,f)<\infty$ and ${D_{\phi}}(f,g)<\infty$ and since $b\in \Theta^{{D_{\phi}}} $, the law of large numbers enables us to state that $|\frac{1}{n}\Sigma_{i=1}^{n}\varphi'\{\frac{f_{a}(\transp aY_i)}{g_{a}(\transp aY_i) }\frac{g(Y_i)}{f(Y_i)}\}\frac{f_{a}(\transp aY_i)}{g_{a}(\transp aY_i) }-\int\ \varphi'(\frac{g(x)}{f(x)}\frac{f_b(\transp bx)}{g_b(\transp bx)})\ g(x)\frac{f_a(\transp ax)}{g_a(\transp ax)}\ dx|\to 0\ a.s.$\\
Furthermore, $|\frac{1}{n}\Sigma_{i=1}^{n}\varphi'\{\frac{f_{a,n}(\transp aY_i)}{g_{a,n}(\transp aY_i) }\frac{g_n(Y_i)}{f_n(Y_i)}\}\frac{f_{a,n}(\transp aY_i)}{g_{a,n}(\transp aY_i) }-\frac{1}{n}\Sigma_{i=1}^{n}\varphi'\{\frac{f_{a}(\transp aY_i)}{g_{a}(\transp aY_i) }\frac{g(Y_i)}{f(Y_i)}\}\frac{f_{a}(\transp aY_i)}{g_{a}(\transp aY_i) }|$

$\ \ \ \ \ \ \ \ \ \ \ \ \ $ $\leq \frac{1}{n}\Sigma_{i=1}^{n}|\varphi'\{\frac{f_{a,n}(\transp aY_i)}{g_{a,n}(\transp aY_i) }\frac{g_n(Y_i)}{f_n(Y_i)}\}\frac{f_{a,n}(\transp aY_i)}{g_{a,n}(\transp aY_i) }-\varphi'\{\frac{f_{a}(\transp aY_i)}{g_{a}(\transp aY_i) }\frac{g(Y_i)}{f(Y_i)}\}\frac{f_{a}(\transp aY_i)}{g_{a}(\transp aY_i) }|$\\
and $|\varphi'\{\frac{f_{a,n}(\transp aY_i)}{g_{a,n}(\transp aY_i) }\frac{g_n(Y_i)}{f_n(Y_i)}\}\frac{f_{a,n}(\transp aY_i)}{g_{a,n}(\transp aY_i) }-\varphi'\{\frac{f_{a}(\transp aY_i)}{g_{a}(\transp aY_i) }\frac{g(Y_i)}{f(Y_i)}\}\frac{f_{a}(\transp aY_i)}{g_{a}(\transp aY_i) }|\to 0$
as a result of the hypotheses initially introduced on $\theta_n.$
Consequently, $\frac{1}{n}\Sigma_{i=1}^{n}|\varphi'\{\frac{f_{a,n}(\transp aY_i)}{g_{a,n}(\transp aY_i) }\frac{g_n(Y_i)}{f_n(Y_i)}\}\frac{f_{a,n}(\transp aY_i)}{g_{a,n}(\transp aY_i) }-\varphi'\{\frac{f_{a}(\transp aY_i)}{g_{a}(\transp aY_i) }\frac{g(Y_i)}{f(Y_i)}\}\frac{f_{a}(\transp aY_i)}{g_{a}(\transp aY_i) }|\to 0$, as it is a  Ces\`aro mean. This enables us to conclude. Similarly, we obtain

$\sup_{b\in\R^d_*}|\frac{1}{n}\Sigma_{i=1}^n\varphi^*\{\varphi'\{\frac{f_{a,n}(\transp aX_i)}{g_{a,n}(\transp aX_i) }\frac{g_n(X_i)}{f_n(X_i)}\}\}-\ \int\ \varphi^*(\varphi'(\frac{g(x)}{f(x)}\frac{f_b(\transp bx)}{g_b(\transp bx)}))f(x)dx|\to0\ a.s.$\hfill$\Box$\\
\noindent{\bf Proof of theorem \ref{LOIDUCRITERE}.}
\begin{theoreme} \label{LOIDUCRITERE} 
Assuming that $(H1)$ to $(H3)$, $(H6)$ and $(H8)$ hold. Then,

$\sqrt n(Var_{\PP}(M(\check c_n(\check \gamma_n),\check \gamma_n)))^{-1/2}(\Pn_nM(\check c_n(\check \gamma_n),\check \gamma_n)-\Pn_nM(a_k,a_k)) \cvL \cN(0,I)$,\\
where $k$ represents the $k^{th}$ step of the  algorithm and with $I$ being the identity matrix in $\R^d$.
\end{theoreme}
\noindent Note that $k$ is fixed in theorem \ref{LOIDUCRITERE} since $\check \gamma_n   =\ arg\inf_{a\in\Theta }\ \sup_{c\in\Theta }\ \Pn_nM(c,a)$ where $M$ is a known function of $k$ , $f$ and $g^{(k-1)}$ - see section \ref{HypoF}.
\begin{proof}
Through a Taylor development of $\Pn_nM(\check c_n(a_k),\check \gamma_n)$ of rank 2, we get at point $(a_k,a_k)$:\\
$\Pn_nM(\check c_n(a_k),\check \gamma_n)$
$=\Pn_nM(a_k,a_k)+\Pn_n\frac{\dr}{\dr a}M(a_k,a_k)\transp {(\check \gamma_n-a_k)}+\Pn_n\frac{\dr}{\dr b}M(a_k,a_k)\transp {(\check c_n(a_k)-a_k)}$

$+\frac{1}{2}\{\transp {(\check \gamma_n-a_k)}\Pn_n\frac{\dr^2}{\dr a\dr a}M(a_k,a_k)(\check \gamma_n-a_k)+\transp {(\check c_n(a_k)-a_k)}\Pn_n\frac{\dr^2}{\dr b\dr a}M(a_k,a_k)(\check \gamma_n-a_k)$

$+\transp {(\check \gamma_n-a_k)}\Pn_n\frac{\dr^2}{\dr a\dr b}M(a_k,a_k)(\check c_n(a_k)-a_k)+\transp {(\check c_n(a_k)-a_k)}\Pn_n\frac{\dr^2}{\dr b\dr b}M(a_k,a_k)(\check c_n(a_k)-a_k)\}$\\
The lemma below enables us to conclude.
\begin{lemme} 
Let $H$ be an integrable function and let $C=\int\ H\ d\PP$ and $C_n=\int\ H\ d\Pn_n$,

$\ \ \ \ \ \ \ \ \ \ \ \ $then, $C_n-C=O_{\PP}(\frac{1}{\sqrt n}).$
\end{lemme}
Thus we get $\Pn_nM(\check c_n(a_k),\check \gamma_n)=\Pn_nM(a_k,a_k)+O_{\PP}(\frac{1}{n}),$
\\ i.e. $\sqrt n(\Pn_nM(\check c_n(a_k),\check \gamma_n)-\PP M(a_k,a_k))=\sqrt n(\Pn_nM(a_k,a_k)-\PP M(a_k,a_k))+o_{\PP}(1).$ \\Hence $\sqrt n(\Pn_nM(\check c_n(a_k),\check \gamma_n)-\PP M(a_k,a_k))$ abides by the same limit distribution as\\ $\sqrt n(\Pn_nM(a_k,a_k)-\PP M(a_k,a_k))$, which is $\cN(0,Var_{\PP}(M(a_k,a_k)))$.\end{proof}
\noindent{\bf Proof of theorem \ref{KernelLOIDUCRITERE}.}
Through proposition \ref{QuotientDonneLoi} and theorem \ref{LOIDUCRITERE}, we derive theorem \ref{KernelLOIDUCRITERE}.\hfill$\Box$\\

\end{document}